\documentclass[11pt,a4paper,leqno]{article}
\usepackage{epsfig}
\usepackage{amssymb}
\usepackage{amsmath}

\author{Hwajeong Kim \\
  \small Humboldt-Universit\"at  zu Berlin\\
        \small     Institut f\"ur Mathematik \\ 
      \small  hjkim@mathematik.hu-berlin.de}
\title{Unstable minimal surfaces of annulus type in manifolds}

\usepackage{latexsym}
\usepackage{epsf,epsfig, xspace} 
\usepackage{amsmath,amsthm,amstext,amsfonts,amssymb,amsxtra,eucal} 
\newtheorem{theorem}{Theorem}[section] 
 
\newtheorem{lemma}{Lemma}[section] 
\newtheorem{remark}{Remark}[section] 
\newtheorem{example}{Example}[section]
\newtheorem{proposition}{Proposition}[section]  
\newcommand{\ar}{A_{\rho}}
\newcommand{\hs}{H^{\frac{1}{2},2}} 
\newcommand{\en}{\enspace} 
\newcommand{\mhr}{\mathcal H _{\rho}}
\newcommand{\hr}{H _{\rho}}
\newcommand{\mh}{\mathcal H}

\newcommand{\dach}{^}
\newcommand{\bqr}{\begin{eqnarray}}
\newcommand{\eqr}{\end{eqnarray}}
\newcommand{\bqrs}{\begin{eqnarray*}}
\newcommand{\eqrs}{\end{eqnarray*}}
\newcommand{\mf}{{\cal F}}
\newcommand{\mfr}{\mf_{\rho}}
\newcommand{\me}{{\cal E}}
\newcommand{\mm}{{\cal M}}

\newcommand{\mj}{{\mathbf J}}

\newcommand{\hab}{h_{\alpha\beta}}
\newcommand{\mon}{\text{mon}}

\newcommand{\osc}{\text{osc}}
\newcommand{\dist}{\text{dist}}
\newcommand{\hfn}{\frac{1}{2},2;0}

\begin{document}

\date{}


\maketitle 

\begin{center}{\large\bf Abstract} \end{center}
 
Unstable minimal surfaces are the unstable stationary points of the Dirichlet integral. In order to obtain unstable solutions, the method of the gradient flow together with the minimax-principle is generally used, an application of which was presented in \cite{s3} for minimal surfaces in Euclidean space. We extend this theory to obtain unstable minimal surfaces in Riemannian manifolds. In particular, we consider minimal surfaces of annulus type.

\section{Introduction}

\footnotetext[1]{This paper is based on my thesis \cite{ki} supervised by Professor Michael Gr\"uter}

For given curves $\Gamma_l\subset N, l=1,...,m$ and $\Gamma:=\Gamma_1\cup\cdots\cup\Gamma_m$, where $(N,h)$ is a Riemannian manifold of dimension $n\ge 2$ with metric $(h_{\alpha\beta})$, we denote the generalized Plateau Problem by $\mathcal P(\Gamma)$.  This deals with minimal surfaces bounded by $\Gamma$, in other words parametrizations $X$ defined on $\Sigma\subset \mathbb R^2$ with $\partial \Sigma=\Gamma$, satisfying the following constraints:
     \begin{enumerate}
          \item[(1)]        
             $\tau_h (X) =0$,
          \item[(2)] 
             $|X_u|_h^2 - |X_v|_h^2 =  \langle X_u, X_v \rangle_h =0$,
          \item[(3)]
             $X|_{\partial \Sigma}$ is weakly monotone and onto $\Gamma$, 
     \end{enumerate} 
where $\tau_{h}:=\Delta X^{\alpha} - \Gamma^{\alpha}_{\beta\gamma}\nabla X^{\beta}X^{\gamma}=0$ is the harmonic equation on $(N,h)$ seen as the Euler-Lagrange equation of the energy functional.\medskip

A regular minimal surface is called unstable if its surface area is not a minimum among neighbouring surfaces with the same boundary.\medskip

Extending the Ljusternik-Schnirelmann theory on convex sets in Banach spaces, a variational approach to unstable minimal surfaces of disc or annulus type in $\mathbb R^n$ was proposed in 1983 (\cite{s1}, see also \cite{s3} \cite{s4}). For the minimal surfaces of higher topological structure in $\mathbb R^n$, it was studied in \cite{js}.\medskip

Recently in \cite{ho}, the existence of unstable minimal surfaces of higher topological structure with one boundary in a nonpositively curved Riemannian manifold was studied by applying the method in \cite{s3}. In particular, the first part of that paper considers the Jacobi field extension operator as the derivative of the harmonic extension.\medskip

In this article, we study unstable minimal surfaces of annulus type in manifolds. The Euclidean case was tackled already in \cite{s4}, and our aim is to generalize this result to manifolds satisfying appropriate conditions. Namely, we will  consider two boundary curves $\Gamma_1, \Gamma_2$ in a Riemannian manifold $(N,h)$ such that one of the following holds.   
   \begin{enumerate}
     \item[(C1)]   
        There exists $p\in N$ with $\Gamma_1, \Gamma_2 \subset B(p,r)$, where  $B(p,r)$ lies within the normal range of all its points. We assume $ r<\pi/(2 \sqrt{\kappa})$, where $\kappa$ is an upper bound of the sectional curvature of $(N,h)$.   
     \item[(C2)]   
        $N$ is compact with nonpositive sectional curvature.     
   \end{enumerate} 
These conditions are related to the existence and uniqueness of the harmonic extension for a given boundary parametrization. \medskip

First, we construct suitable spaces of functions, the boundary parametrizations, distinguishing the cases (C1) and (C2). We introduce a convex set which serves as a tangent space for the given boundary parametrization. Then we consider  the following functional: 
    \[\me(x):= \frac{1}{2}\int |d\mf(x)|_h^2,\]
where $\mf(x)$ denotes the harmonic extension of annulus type or of two-disc type with boundary parametrization $x$. We next discuss the differentiability of $\me$, in particular for the case in which the topology of the surfaces changes (from an annulus to two discs). Defining critical points of $\me$, will show the equivalence between the harmonic extensions (in $N$) of critical points of $\me$ and minimal surfaces in $N$. The $H^{2,2}$-regularity of the harmonic extension of a critical point of $\me$ (see the appendix or \cite{ki2}) plays an important role in the argument.\medskip

In section \ref{CPT}, we prove the Palais-Smale condition for $\me$. In particular, we investigate carefully the behaviour of boundary mappings which are fixed at only one point. In order to deform level sets of $\me$, we also construct a suitable vector field and its corresponding flow. \medskip

Roughly speaking, Lemma \ref{boundary} shows that the energy of some annulus-type harmonic extensions is greater than that of two-disc type harmonic extensions by a uniformly positive constant. Although this result refers to Riemannian manifolds, it turns out to be more restrictive than that of Euclidean spaces, which holds uniformly on any bounded set of boundary parametrizations. This somewhat weaker result is anyhow enough for the present purposes.\medskip 

Following the arguments set out in \cite{s1}, we can prove the main theorem of this paper. This states that if there exists a minimal surface (of annulus type) whose energy is a strict relative minimum in $\mathcal S(\Gamma_1,\Gamma_2)$ (suitably defined for each case (C1) and (C2)), the existence of an unstable minimal surface of annulus type is ensured under certain assumptions related to the solutions of $\mathcal P(\Gamma_i)$. We eventually apply this result to the three-dimensional sphere $S^3$ and the three-dimensional hyperbolic space $H^3$, whose curvatures are $1$ and $-1$, respectively.


\section{Preliminaries} 


\subsection{Some definitions}

Let $(N, h)$ be a connected, oriented, complete Riemannian manifold of dimension $n \ge 2$, embedded isometrically and properly into some $\mathbb R^k$ as a closed submanifold by means of the map $\eta$ (\cite{gr}). Moreover, $d\omega$ and $d_0$ denote the area elements in $\Omega \subset \mathbb R^2$ and in $\partial \Omega$ respectively.\medskip

Indicating  
    \[B:= \{w \in \mathbb R^2 \mid |w| < 1 \}\]
we define  
         \[ H^{1,2}\cap C^0(B,N) := \{ f\in H^{1,2}\cap C^0(B,\mathbb R^k) | f(B)\subset N \} \]
with norm $ \| f \|_{1,2 ;0}:= \|d f\|_{L^2} + \|f\|_{C^0}$.  Now set  
        \[T_f H^{1,2}\cap C^0(B,N) \cong \{ V\in H^{1,2}\cap C^0(B,\mathbb R^k) | V(\cdot) \in T_{f(\cdot)} N \} 
                                =: H^{1,2}\cap C^0(B,f^{\ast}TN), \]
with norm 
  \bqr
      \|V\|  :=  \big( \int_B |\nabla^f V|^2_h d\omega \big)^{\frac{1}{2}} + \|V\|_{C^0} 
             \cong   \big( \int_B |dV|^2_{\mathbb R^k} d\omega \big)^{\frac{1}{2}} + \|V\|_{C^0}.
           \label{yeoseot}
     \eqr 

Let  $\Gamma$ be a Jordan curve in $N$ diffeomorphic to $S^1:=\partial B$. Then $N$ can be equipped with another metric $\tilde h$ such that $\Gamma$  is a geodesic in $(N, \tilde{h})$. We observe that $H^{1,2}\cap C^0\big((B,\partial B),(N,\Gamma)_{\tilde h} \big)$ and $H^{1,2}\cap C^0\big((B,\partial B),(N,\Gamma)_{h} \big)$ coincide as sets.\medskip
 
Using the exponential map in $(N,\tilde{h})$, we let    
       \bqrs
          \hs\cap C^0(\partial B;\Gamma) & := & 
          \{ u\in \hs \cap C^0(\partial B, \mathbb R^k) | u(\partial B) = \Gamma \},
       \eqrs
where the norm is given by $ \| u \|_{\frac{1}{2},2 ;0}:= \|d\mh(u)\|_{L^2} + \|u\|_{C^0}$, and $\mh(u)$ is the harmonic extension in $\mathbb R^k$ with $\mh(u)|_{\partial B}(\cdot)=u(\cdot)$. In addition 
       \bqrs
           T_u \hs\cap C^0(\partial B;\Gamma) & := & \{ \xi \in  \hs\cap C^0(\partial B, u^{\ast}TN) | 
                                                        \xi(z) \in T_{u(z)}\Gamma, 
                                                \en  \text{for all} \en z\in \partial B\} \\
                                               & = & \hs\cap C^0(\partial B,u^{\ast}T\Gamma). 
       \eqrs
Finally, the energy of $f\in H^{1,2}(\Omega, N)$ is denoted by 
    \[E(f) := \frac{1}{2}\int_{\Omega}|df|_h^2dw.\]


\subsection{The setting}\label{construction}

Let $\Gamma _1, \Gamma_2$ be two Jordan curves of class $C^3$ in $N$ with diffeomorphisms $\gamma^i : \partial B \rightarrow \Gamma_i, i=1,2$, and $\dist(\Gamma _1, \Gamma_2)>0$. For $\rho \in (0,1)$ let
    \[ A_{\rho} := \{w \in B \mid \rho < |w| < 1 \}\]
have boundary $C_1 := \partial B$ and $C_\rho  := \partial B_{\rho} =:C_2 \,(\rho \en\text{fixed})$, and indicate    
    \[\mathcal X^i_{\mon}:=\{x^i\in \hs\cap C^0(\partial B ; \Gamma_i)\,|\,x^i\en\text{is 
                  weakly monotone and onto} \en \Gamma_i \en \text{with degree}\en 1\}.\]

{\bf I)} We first consider the following condition for $(N,h)(\supset \Gamma_1, \Gamma_2)$.
   \begin{enumerate}
     \item[(C1)] 
           There exists $p\in N$ with $\Gamma_1, \Gamma_2 \subset B(p,r)$, where  $B(p,r)$ lies within the normal range of all its points. We assume $ r<\pi/(2 \sqrt{\kappa})$, where $\kappa$ is an upper bound of the sectional curvature of $(N,h)$.   
   \end{enumerate} 

Throughout the paper, $B(p,r)$ denotes a geodesic ball with center $p\in N$ as in (C1).\medskip
 
We can easily observe the following property (see \cite{ki2}).

\begin{remark}\label{trace}
 If $\Gamma_1, \Gamma_2\subset N$ satisfy (C1), then for each $x^i\in \hs \cap C^0(\partial B;\Gamma_i)$ and $\rho\in(0,1)$ there exist  $g_{\rho}\in H^{1,2}\cap C^0(\overline{\ar}, B(p,r))$ and $g^i\in H^{1,2}\cap C^0(\overline{B}, B(p,r))$ with $g_{\rho}|_{C_1} = x^1,\, g_{\rho}|_{C_{\rho}}(\cdot) = x^2(\frac{\cdot}{\rho})$ and $g^i|_{\partial B}=x^i, \, i= 1,2$.
\end{remark}

 From the results in \cite{hkw}, \cite{jk} and the above remark, we have a unique harmonic map of annulus and disc type in $B(p,r)\subset N$ for a given boundary mapping in the class $\hs\cap C^0$. Now we define  
        \bqrs 
           M^i 
               & := & \{x^i\in \mathcal X^i_{\mon}\,|\,x^i \en\text{preserves the orientation} \}. 
        \eqrs
Then $M^i$ is complete, since the $C^0$-norm preserves monotonicity.\medskip

Moreover, let
     \bqrs   
        \mathcal S (\Gamma_1,\Gamma_2) & = &  \{ X\in H^{1,2}\cap C^0(\overline{A_{\rho}},B(p,r)) |
                                                 \, 0<\rho<1,  \en X|_{C_i} \text{is weakly monotone}\},\\
        \mathcal S (\Gamma_i) & = &  \{ X\in H^{1,2}\cap C^0(\overline{B},B(p,r))| X|_{\partial B}
                                                    \en \text{is weakly monotone}\}.    
     \eqrs

{\bf II)} We now investigate another significant condition for $(N,h)$.
 \begin{enumerate}
     \item[(C2)]   
          $N$ is compact with nonpositive sectional curvature.     
  \end{enumerate} 
 
A compact Riemannian manifold is homogeneously regular and the condition of nonpositive sectional curvature implies $\pi_2(N) = 0$.\medskip

In order to define $M^i$, we first consider for $\rho\in(0,1)$ the following 
 \begin{equation*} 
     \tilde{G_ {\rho}} := \{f \in H^{1,2}\cap C^0(\overline{\ar},N) |\, f|_{C_i}\, \text{is continuous}, \en \text{weakly monotone and onto} \en \Gamma_i\}.
 \end{equation*}

We may take a continuous homotopy class, denoted by $\tilde{F_ {\rho}}\subset\tilde{G_ {\rho}}$, so that every two elements $f, g$ in $\tilde{F_ {\rho}}$ are continuously homotopic $f\sim g$ (not necessarily fixing the boundary parametrization). We further demand some relation $\tilde{F_ {\rho}} \sim \tilde{F_ {\sigma}}$ to hold for any $\rho, \sigma \in(0,1)$. Precisely, for some $\tilde f \in \tilde{F_ {\sigma}}$, $f \in \tilde{F_ {\rho}}$ and some diffeomorphism  $\tau^{\rho}_{\sigma}: [\sigma, 1] \rightarrow [\rho, 1]$, we require $\tilde f(r,\theta) = f(\tau^{\rho}_{\sigma}(r),\theta)$. Let $\tilde{F_ {\rho}}$ be fixed. Then for any $\sigma\in (0,1)$ we can find $\tilde{F_ {\sigma}}$ with  $\tilde{F_ {\rho}} \sim \tilde{F_ {\sigma}}$.\medskip

We now consider all possible $H^{1,2}\cap C^0$-extensions of disc type in $N$: 
  \[ \mathcal S (\Gamma_i)  :=  \{ X\in H^{1,2}\cap C^0(\overline{B},N) | X|_{\partial B}\en 
     \text{is weakly monotone  onto}\en \Gamma_i\},\]
assuming that this set is not empty, for each $i=1,2$.

\begin{lemma} \label{hundredone}
   \begin{enumerate}
     \renewcommand{\labelenumi}{(\roman{enumi})}
     \item
       For $X^1\in \mathcal S (\Gamma_1)$ and $X^2\in \mathcal S (\Gamma_2)$ 
         there exists $f_{\rho}\in H^{1,2}\cap C^0(\ar,N)$ such that 
         $f_{\rho}|_{C_1}(\cdot) = X^1|_{\partial B}(\cdot)$ and 
         $f_{\rho}|_{C_\rho}(\cdot) = X^2|_{\partial B}(\frac{\cdot}{\rho})$, 
         for $\rho\in(0,1)$. 
     \item
      Moreover, there exists $\rho_0\in(0,1)$ and a uniform positive constant $C$ such that for some 
        $f_{\rho}\in H^{1,2}\cap C^0(\ar,N)$, with 
        $f_{\rho}|_{C_\rho}(\cdot) = X^2|_{\partial B}(\frac{\cdot}{\rho})$  
          \bqr\label{hundredthree}
            E(f_{\rho})\le C, \en \text{for all}\en \rho\le \rho_0.
          \eqr 
   \end{enumerate}    
\end{lemma} 
{\bf Proof.} {\it (i)} For a given $\varepsilon >0$, take $\sigma_i>0$ with $\osc_{B_{\sigma_i}} X^i < \varepsilon$. Choose $\rho>0$ with $\frac{\rho}{\sigma_2} < \sigma_1$, and let $\mathcal H : B_{\sigma_1}\backslash B_{\frac{\rho}{\sigma_2}} \rightarrow \mathbb R^k$ be harmonic with $X^1|_{\partial B_{\sigma_1}} - X^1(0)$ on $\partial B_{\sigma_1}$ and $X^2|_{\partial B_{\sigma_2}} - X^2(0)$ on $\partial B_{\frac{\rho}{\sigma_2}}$. This implies $\| \mathcal H \|_{C^0} < \varepsilon $. Now let $g\in H^{1,2}\cap C^0 (B_{\sigma_1}\backslash B_{\frac{\rho}{\sigma_2}}, N)$ with $X^1(0)$ on $ \partial B_{\sigma_1} $ and $X^2(0)$ on $ \partial B_{\frac{\rho}{\sigma_2}}$. \medskip

Considering coordinate neighbourhoods for the submanifold $N\stackrel{\eta}{\hookrightarrow} \mathbb R^k$, we may take a finite covering of $f_{\rho}((\overline{\ar}))$, and by projection we obtain a smooth map $r:\mathcal N_{\delta}(f_{\rho}(\overline{\ar}))\rightarrow N$ with $r|_{\mathcal N_{\delta}(f_{\rho}(\overline{\ar}))\cap N} = Id$ for some $\delta>0$, where $\mathcal N_{\delta}(\cdot)$ is $\delta$-neighbourhood in $\mathbb R^k$. Setting $T(s,\theta):= (\frac{1}{s}\rho, \theta)$ in polar coordinates, we can define $f_{\rho}$ with the desired properties: 
          \bqr\label{hundred}
                                            f_{\rho}:= \left\{ \begin{array}{r@{\quad,\quad} l}
                     X^1|_{B\backslash B_{\sigma_1}} & \text{on}\en B\backslash B_{\sigma_1},\\
                               r\circ(g+ \mathcal H) & \text{on}\en  B_{\sigma_1}\backslash 
                                                            B_{\frac{\rho}{\sigma_2}}, \\
                                  X^2(T^{-1}(\cdot)) & \text{on}\en B_{\frac{\rho}{\sigma_2}}
                                                            \backslash B_{\rho}.
                                                    \end{array} \right. 
          \eqr           
      
   {\it(ii)} The claim follows from the above construction, since $\frac{\rho}{\sigma_2} < \sigma_1, \,\rho\le \rho_0$ for some $\rho_0>0$. \hfill $\Box$ \medskip   
       
Under the assumption that $\mathcal S(\Gamma_i)\not= \emptyset$, for given $\Gamma_i\in N$ we have an annulus-type-extension like that of (\ref{hundred}), and we take homotopy classes which contain such an extension. From now on twiddles will be dropped. \medskip 

Define
        \bqr
           \mathcal S(\Gamma_1, \Gamma_2) &:=&  \{f\in F_{\rho}\,|\, 0 < \rho < 1\},\label{mathcals}
        \eqr
as well as the two function spaces
      \begin{eqnarray*}
       M^1 &:=&\{x^1(\cdot)=f|_{C_1}(\cdot ),
         \,f\in \mathcal S(\Gamma_1, \Gamma_2) |\,x^1\en\text{is orientation preserving with degree 1}\},\\
       M^2 &:=& \{x^2(\cdot)=f|_{C_{\rho}}(\cdot \rho),
         \,f\in \mathcal S(\Gamma_1, \Gamma_2) |\,x^2\en\text{is orientation preserving with degree 1}\}.
      \end{eqnarray*}
   
For $x^i \in \mathcal X^i_{\mon}$, $\mathcal H_{\rho}(x^1,x^2)$ denotes the unique $\mathbb R^k$-harmonic extension on $A_\rho$ with $x^1(\cdot)$ on $C_1$ and $x^2(\frac{\cdot}{\rho})$ on $C_{\rho}$, while $\mathcal H(x)$ is the $\mathbb R^k$-harmonic extension of disc type with boundary $x\in \mathcal X^i_{\mon}$.

\begin{lemma}\label{tangent}
   \begin{enumerate}
     \renewcommand{\labelenumi}{(\roman{enumi})}
        \item \label{tangent2}
          For each $x^i_0 \in  {M^i}$, $i=1,2$, there exists $\varepsilon(x^i_0)>0$ such that  
         \[\text{if} \en x^i\in \mathcal X^i_{\mon}\en\text{with}\en
             \|x^i - x^i_0\|_{\frac{1}{2},2;0} < \varepsilon, \en \text{then}\en x^i\in  {M^i}. \]
        \item 
             $ {M^i}$ is complete with respect to $\|\cdot\|_{\frac{1}{2},2;0}$. 
        \end{enumerate}
 \end{lemma}
{\bf Proof.} {\it(i)} Let $f_\rho \in \tilde{F_\rho}$ with $f_\rho|_{C_1} = x^1_0$ and $ f_\rho|_{C_{\rho}}(\cdot) = y^2(\frac{\cdot}{\rho})$ for some $y^2 \in  {M^2} $.\medskip
 
We consider the smooth retraction $r:\mathcal N_{\delta}(f_{\rho}(\overline{\ar}))\rightarrow N$ as in the proof of Lemma \ref{hundredone}. Let $\|x^i - x^i_0\|_{\frac{1}{2},2;0} < \varepsilon<\delta$. Then by Lemma 4.2 from \cite {s4}, 
       \bqrs            
            \lefteqn { \int _{\ar}|d(r(f_\rho + \mhr(x^1 - x^1_0, 0)))|^2d\omega } \\
               & & \le C(\|f_\rho\|_{C^0},\varepsilon,N)\big(\int_{\ar}|df_{\rho}|^2d\omega + 
                          \int_B|d\mathcal H(x^1-x^1_0)|^2d\omega\big)
                        \le C(\|f_\rho\|_{1,2;0},\varepsilon,N).
         \eqrs
       
Now, let $H(t,\cdot):= (1-t)\mhr(x^1-x^1_0,0): [0,1]\times \ar \rightarrow \mathbb R^k$ with $\|H\|_{C^0} < \varepsilon$ and  $G :[0,1]\times \ar \rightarrow N$ with  $G(t,\cdot)= f_{\rho}(\cdot)$ for all $t\in  [0,1]$. Since $r(G+H) : [0,1]\times \ar \rightarrow N$ is a homotopy between $f_{\rho}$ and $r(f_{\rho}+\mhr(x^1 - x^1_0, 0))$, it follows $r(f_\rho + \mhr(x^1 - x^1_0, 0))(\sim f_{\rho})\in \tilde{F_\rho}$, and $x^1 \in  {M^1}$. Similarly, we can prove that $x^2 \in {M^2}$ if $\|x^2 - x^2_0\|_{\frac{1}{2},2;0}  < \varepsilon'$ for some small $\varepsilon'>0$.\medskip

{\it(ii)} A Cauchy sequence $\{x^i_n\}\subset M^i$ converges to $x^i\in\hs \cap C^0 (\partial B;\Gamma_i)$, and for some $n$, $\|x^i_n- x^i\|_{C^0}< \varepsilon$. Considering  $\mhr(x^1 - x^1_n, 0)$ and $g_\rho\in F_{\rho}$ with $x^1_n$ on $C_1$ and $0$ on $C_{\rho}$, we can find a homotopy in $N$ between $g_{\rho}$ and $r(g_{\rho}+\mhr(x^1 - x^1_0, 0)) $ as in {\it(i)}. We may also apply this argument to $x^2$. Note that $x^i$ is weakly monotone, and hence $x^i\in M^i$, concluding the proof. \hfill $\Box$ \medskip

From the proof we easily conclude that the set of $x^i$'s which possess annulus-type-extensions with uniform energy with respect to $\rho\le \rho_0$ is an open and closed subset of $\mathcal X^i_{\mon}$. Thus, it is a non-empty connected component of $\mathcal X^i_{\mon}$ and must coincide with $M^i$, since $M^i$ is a connected subset of $\mathcal X^i_{\mon}$. Hence we obtain the following property.
   
\begin{remark}\label{hundredtwo}
 For each $x^i\in M^i, i=1,2$, there exist $f_{\rho}\in \mathcal S(\Gamma_1, \Gamma_2)$ and $C>0$ with $E(f_{\rho})\le C$ for all $\rho\le\rho_0$ and some $\rho_0\in(0,1)$. Clearly, this result also holds for $x^i\in M^i$ if we assume (C1).
\end{remark}

For disc-type extensions of $x^i\in M^i$ the following lemmata will be useful.

 \begin{lemma} \label{mo}
Let $(N,h)$ be a homogeneously regular manifold and $u$ an absolutely continuous map of $\partial B_{r}(x_0)$ into $N\ni x_0$ with $\int^{2\pi}_{0}|u'(\theta)|^2_{h}d\theta \le \frac{C'}{\pi}$. Then there exists $f\in H^{1,2}(B_{r}(x_0),N)\cap C^0(\overline{B_{r}(x_0)},N)$ with $f|_{\partial B_{r}(x_0)} = u$ and $ E_{B_{r}(x_0)}(f) \le \frac{C''}{C'}\int^{2\pi}_{0}|u'(\theta)|^2_{h}d\theta$, where $C'', C'$ are the constants defined by homogeneous regularity.
\end{lemma}
{\bf Proof.} See \cite {mo} Lemma 9.4.8 b).                                              \hfill $\Box$ 

\begin{lemma} \label{cou.leb}
Let $f_{\rho} \in H^{1,2}(\ar,N)$, $0<\rho<1$. For each $\delta \in (\rho,1)$ there exists $\tau \in (\delta, \sqrt{\delta})$ with $\int^{2\pi}_{0}\left|\frac{\partial f_{\rho}(\tau,\theta)}{\partial \theta}\right|^2_hd\theta \le \frac{4E(f_{\rho})}{\ln \frac{1}{\delta}}$. 
\end{lemma}
 {\bf Proof.} Similar to the proof of the Courant-Lebesgue Lemma.       \hfill $\Box$\medskip

For $x^i\in M^i$, and given the choice of $\mathcal S(\Gamma_1,\Gamma_2)$, Remark \ref{hundredtwo} tells that we can find $f_{\rho}\in H^{1,2}(\ar,N)$ with boundary $x^i$ such that $E(f_{\rho})\le C$ for all $\rho\le\rho_0$. Then from Lemma \ref{cou.leb} and Lemma \ref{mo}, we have  $g_{\tau}\in H^{1,2}(B_{\tau},N)$ with boundary $f_{\rho}|_{\partial B_{\tau}}$ for some $\rho$. Together with $g_{\tau}$ and $f_{\rho}|_{B\backslash B_{\tau}}$, we obtain a map $X \in H^{1,2}(B,N)$ with boundary $x^1$. Similarly, we have $\tilde X\in H^{1,2}(B,N)$ with boundary $x^2$.\medskip

Moreover, the harmonic extension of disc type for each $x^i\in M^i$ in $N$ is unique, independently of the choice of homotopy class $\mathcal S(\Gamma_1, \Gamma_2)$, because of the following well known fact.

\begin{lemma}\label{disc}
    $\pi_2(N)=0 \Leftrightarrow$ Any $h_0, h_1 \in C^0(B,N)$ with $h_0 |_{\partial B} = h_1 |_{\partial B}$ 
    are homotopic. 
\end{lemma}

On the other hand, using the construction (\ref{hundred}) and the previous Lemma we can easily check that the traces of the elements in $\mathcal S(\Gamma_i)$ belong to $M^i$. From \cite {es}, \cite{le}, \cite{hm}, we then have the following. 
   
\begin{remark}
   \begin{enumerate}
     \renewcommand{\labelenumi}{(\roman{enumi})}
       \item 
          For $x^i \in M^i$, there exist a unique harmonic extension of disc type on $B$ and of annulus type on $\ar$, $\rho\in (0,1)$.
       \item
         The elements of $M^i$ are the traces of the elements of $\mathcal S(\Gamma_i)$.
    \end{enumerate}
 \end{remark}

{\bf III) \,Now let $(N,h)$ and $\Gamma_i, i=1,2$ satisfy (C1) or (C2).}\medskip

Observing $\partial B \cong \mathbb R/2\pi$, for a given oriented $y^i\in\mathcal X^i_{\mon}$ there exists a weakly monotone map $w^i\in C^0(\mathbb R, \mathbb R)$ with $w^i(\theta + 2\pi) = w^i(\theta) + 2\pi$ such that $y^i(\theta) = \gamma^i (\cos(w^i(\theta)), \sin(w^i(\theta)))=:\gamma^i\circ w^i(\theta)$. In addition $w^i = \tilde{w}^i + Id$ for some $\tilde{w}^i\in C^0(\partial B,\mathbb R)$.\medskip 
    
Denoting the Dirichlet integral by $D$ and the  $\mathbb R^k$-harmonic extension by $\mathcal H$, let 
   \begin{equation*}
     W^i_{\mathbb R^k} := \{ w^i\in C^0(\mathbb R, \mathbb R)\,|\,w^i\en \text{is weakly monotone},    
               w^i(\theta + 2\pi) = w^i(\theta) + 2\pi ; 
               D(\mathcal H (\gamma^i\circ w^i)) < \infty \}. 
    \end{equation*}
 Clearly, $W^i_{\mathbb R^k}$ is convex (for further details, refer to \cite{s1}). \medskip
  
 Now take $x^i\in M^i$. Considering $w - w^i$ as a tangent vector along $\tilde{w}^i$, let     
    \bqrs 
     \mathcal T_{x^i} = \{ d\gamma^i((w - w^i)\frac{d}{d\theta}\circ\tilde{w}^i)
            \,|\,w\in W^i_{\mathbb R^k}\en\text{and}\en \gamma^i\circ w\dach i =x\dach i\}.
   \eqrs
Note that $\mathcal T_{x^i}$ is convex in $T_{x^i}\hs\cap C^0(\partial B;\Gamma_i)$, since $W^i_{\mathbb R^k}$ is convex. For $\xi=d\gamma^i((w - w^i)\frac{d}{d\theta}\circ\tilde{w}^i)\in\mathcal T_{x^i}$ we have that $\widetilde{\exp}_{x^i}\xi = \gamma^i(w)$, $\widetilde{\exp}$ denoting the exponential map with respect to the metric $\tilde h$.  \medskip
 
If (C1) holds, then clearly $\widetilde{\exp}_{x^i}\xi\in M^i$ for $\xi\in \mathcal T_{x^i}$. For the case (C2), let us recall the proof of Lemma \ref{tangent}. Since $N$ is compact, there exists $l_i>0$, depending on $\gamma^i$, such that for any $x^i\in M^i$, $\widetilde{\exp}_{x^i}\xi \in M^i,\en \text{provided that}\en \|\xi\|_{\mathcal T_{x^i}} < l_i$.\medskip

The following set-up holds true in both cases (C1) and (C2).\medskip

\newpage

 {\bf Definition} 
    \begin{enumerate}
        \renewcommand{\labelenumi}{(\roman{enumi})}
          \item
            Let $\mathcal M := M^1 \times M^2 \times (0,1)$ with the product topology and $x := (x^1,x^2,\rho)\in \mathcal M$. Then the set $\mathcal T_{x}\mathcal M := \mathcal T_{x^1}\times \mathcal T_{x^2} \times \mathbb R$ is convex.\smallskip
       
Let $\mf(x) = \mf(x^1,x^2,\rho) = \mf_{\rho}(x^1,x^2):\ar \rightarrow N$ be the unique harmonic extension with $x^1$ on $C^1$ and $x^2(\frac{\cdot}{\rho})$ on $C^2$, and define 
            \bqrs
              \me &:& \mm \longrightarrow \mathbb R\\ 
                  & & x \longmapsto 
                     E(\mf(x)):=  \frac{1}{2}\int_{\ar}|d\mf_{\rho}(x^1,x^2)|_h^2d\omega.
            \eqrs
         \item
Define $\partial\mm := M^1\times M^2\times \{0 \}$, $ \mathcal T_x\partial\mm := \mathcal T_{x^1}\times \mathcal T_{x^2}$ and $ \overline{\mm} := \mm \cup \partial\mm$.\smallskip
  
Let $\mf^i(x^i):\ar \rightarrow N$ be the unique harmonic extension with boundary $x^i$, for $x = (x^1,x^2,0) \in \partial \mm$, and define
            \[ \me(x)  :=  E(\mf^1(x^1)) + E(\mf^2(x^2)).\]
    \end{enumerate}


\subsection{Harmonic extension operators}

Let $\Omega=\ar$ or $\Omega=B$.  A weak Jacobi field $\mj$ with boundary $\xi$ along a harmonic function $f$ is a weak solution of  
   \[ \int_{\Omega} \langle\nabla \mj, \nabla X \rangle + \langle tr\, R(\mj, df)df, X \rangle d\omega = 0,\] 
for all X $\in H^{1,2}(\Omega, f^{\ast}TN)$ with $X|_{\partial \Omega}=\xi$. Weak Jacobi fields are natural candidate derivatives of the harmonic operators  $\mfr$ and $\mf^{i}$.\medskip

We have the following property of weak Jacobi fields, from \cite{ho}.

\begin{lemma}\label{jacobifield} 
The weak Jacobi field $\mj$ with boundary $\eta\in T_{x^i}\hs\cap C^0$ along a harmonic $\mf$ with boundary $x^i$ is well defined in the class $H\dach{1,2}$ and continuous up to the boundary. It satisfies  
      \bqrs
       \| \mj_{\mf}\|_{C^0} \le  \|\mj_{\mf}|_{\partial \Omega}\|_{C^0}, \en 
       \| \mj_{\mf}\|_{1,2;0 }  
                  \le   C(N, \|f\|_{1,2:0})\|\mj_{\mf}|_{\partial \Omega}\|_{\frac{1}{2},2 ;0 }.    
     \eqrs
 \end{lemma}

Now we can discuss the differentiability of harmonic extension operators. 

\begin{lemma} \label{harmonic} 
The operators $\mfr, \mf^i$ are partially differentiable in $x^1$ (resp. $x^2$) for variations in $T_{x^1}\hs\cap C^0$ (resp. $T_{x^2}\hs\cap C^0$). Their derivatives are continuous Jacobi field operators with respect to $x^1, x^2$.
\end{lemma}
{\bf Proof.} The proof reproduces an argument we shall explain in full detail in Lemma \ref{differ}, cases (B), (C), and as such will not be anticipated here. Alternatively, one can follow the aforementioned \cite{ho}. \hfill $\Box$


\section{The variational problem}\label{differentiabilitycritical}  
\subsection{Differentiability of $\me$ on $\overline{\mm}$}\label{defferentiability}

\begin{lemma}\label{differ}
  The following hold: 
   \begin{enumerate}
     \renewcommand{\labelenumi}{(\Alph{enumi})}
       \item
          $\me$ is continuously partially differentiable in $x^1, x^2$ with respect to variations in $\mathcal T_{x^1}$, $\mathcal T_{x^2}$ and the derivatives are continuous on $M^1\times M^2$. 
       \item 
          $\me$ is continuous with respect to $\rho \in [0,1)$, even uniformly on $\mathcal N_{\varepsilon}(x^i_0)$ for some $\varepsilon>0$ independent of $x^i_0\in M^i$, $i=1,2$.
       \item 
          The partial derivatives in $x^1, x^2$ are continuous with respect to $\rho \in [0,1)$, uniformly continuous on $\mathcal N_{\varepsilon}(x^i_0)$ for some $\varepsilon>0$ independent of $x^i_0\in M^i$, $i=1,2$.
       \item
          $\me$ is differentiable with respect to $\rho\in (0,1)$. 

  \end{enumerate}
\end{lemma}
{\bf Proof.} From now on, continuity will be understood in the sense of convergence of subsequences.\medskip
    
{\it (A)} The Dirichlet integral functional is in $C^{\infty}$, so Lemma \ref{harmonic} guarantees that $\me$ is continuously partially differentiable with continuous partial derivatives on $M^1\times M^2$.\medskip
                             
\underline{Computation of the derivatives}:\smallskip

Let $x=(x^1,x^2,\rho) \in \mm$, $\xi^1 \in \mathcal T_{x^1}$. By Lemma \ref{tangent} there is a small $t_0>0$ such that $\widetilde{\exp}_{x^1}(t\xi^1)\in M^1,\en 0\le t\le t_0$. Thus, 
          \bqr
                 \langle \delta_{x^1}\me,\xi^1\rangle  
                   & := &\frac{d}{dt}\Big|_{t=0}\me(\widetilde{\exp}_{x^1}(t\xi^1),x^2,\rho)\nonumber \\
                   & = &\int_{\ar}\langle d\mfr(x^1,x^2),\nabla D_{x^1}\mfr(x^1,x^2)(\xi^1)\rangle _hd\omega
                                                             \nonumber \\
                   & = &\int_{\ar}\langle d\mfr(x^1,x^2),\nabla\mj_{\mfr}(\xi^1,0)\rangle _hd\omega  
                                                 \hspace{1.0cm}        \text{(by Lemma \ref{harmonic})},
          \eqr
since by computation we obtain, with $\mfr(t):=\mfr\left(\widetilde{\exp}_{x^1}(t\xi^1),x^2\right)$,
               \bqrs
                  \nabla_{\frac{d}{dt}}\left(\mf^{\alpha}_{\rho,i}(t)dx^i\otimes
                                              \frac{\partial}{\partial y^{\alpha}}\circ\mfr(t)\right)
                   = \nabla\frac{d}{dt}\mfr\left(\widetilde{\exp}_{x^1}(t\xi^1),x^2\right)
                 ( =  \nabla\left(D_{x^1}\mfr(x^1,x^2)(\xi^1)\right), t = 0).
               \eqrs 
For $\xi^2 \in \mathcal T_{x^2}$ Lemma \ref{harmonic} yields $\langle \delta_{x^2}\me,\xi^2\rangle = \int_{\ar}\langle d\mfr(x^1,x^2),\nabla\mj_{\mfr}(0,\xi^2(\frac{\cdot}{\rho}))\rangle_hd\omega$. Similarly, for $x = (x^1,x^2,0)\in \partial \mm$, $\langle \delta_{x^i}\me,\xi^i\rangle = \int_B\langle d\mf^i(x^i),\nabla\mj_{\mf^i}(\xi^i)\rangle _hd\omega, i=1,2$.\bigskip 
                           
 
For (B) we shall split the proof into three sub-steps {\it B}-I), {\it B}-II), {\it B}-III). Similarly for (C) we shall have  {\it C}-I), {\it C}-II), {\it C}-III).\medskip

{\it B}-I) The set-up.\medskip

The claim is that $\me$ is continuous when $\rho\rightarrow \rho_0$. Fixing $\rho_0=0$ is no great restriction, since the proof for $\rho_0\in(0,1)$ carries over in an analogous, even easier, fashion. Taking $\rho_0=0$ translates our claim into  
          \begin{equation}\label{hana}
              \int_{\ar}\left|d\mfr(x^1,x^2)\right|^2_hd\omega 
                    \longrightarrow \int_B\left|d\mf^1(x^1)\right|^2_hd\omega
                               +\int_B\left|d\mf^2(x^2)\right|^2_hd\omega
          \end{equation}
uniformly on $\mathcal N_{\varepsilon}(x^i_0)$ for some $\varepsilon>0$ independent of $x^i_0\in M^i$, whenever $\rho\rightarrow 0$.\medskip
      
Let $\mfr := \mfr(x^1,x^2)$ and $\mf^i := \mf^i(x^i), i=1,2$. By Lemma \ref{cou.leb}, for each $\delta $ with $0<\rho<\delta<1$ there exists $\nu \in (\delta,\sqrt{\delta})$ such that 
          \bqr 
           \int^{2\pi}_{0}\left|\frac{\partial \mfr(\nu,\theta)}{\partial \theta}\right|_hd\theta \le
               \sqrt{2\pi}\left(\int^{2\pi}_{0}\left|\frac{\partial \mfr(\nu,\theta)}{\partial\theta}
                                                \right|^2_h d\theta\right) 
                  ^{\frac{1}{2}}\le \frac{C}{\sqrt{|\ln\delta|}}.\label{circle}
          \eqr 
Due to Remark \ref{hundredtwo}, $C$ is independent of $\rho\le\rho_0$, for some $\rho_0\in (0,1)$.\medskip
 
By means of $\mfr$ we now construct two maps by setting  
         \bqr
            f_{\nu} : A_{\nu} \longrightarrow N \qquad &\text{with}& \quad f_{\nu}(re^{i\theta})
                                :=\mfr(re^{i\theta}),
                                           \en re^{i\theta} \in A_{\nu},\nonumber\\
            g_{\nu'} : A_{\nu'} \longrightarrow N \qquad
              &\text{with}&  \quad g_{\nu'}(re^{i\theta}) := \mfr(T(re^{i\theta})),
                                                 \en   re^{i\theta} \in A_{\nu'}.\label{g}
         \eqr  
The constants $\nu':=\frac{\rho}{\nu}$, $\nu\in(\delta,\sqrt{\delta})$ and  $\delta \in (\rho,1)$ satisfy the property (\ref{circle}) in the limit $\nu', \nu\rightarrow 0$ for $\rho \rightarrow 0$. (One can take for instance $\delta = \sqrt{\rho}$). The map $T(re^{i\theta}) = \frac{\rho}{r}e^{i\theta}$ goes from $ A_{\nu'}$ to $B_{\nu}\backslash B_{\rho}$ surjectively. Then, $f_{\nu}$ and $g_{\nu'}$ are harmonic maps into $N$ with $f_{\nu}|_{\partial B} = x^1$, $g_{\nu'}|_{\partial B} = x^2$ and $\osc_{\partial B_{\nu}} f_{\nu}\rightarrow 0,\en \osc_{\partial B_{\nu'}} g_{\nu'} \rightarrow 0$ as $\rho \rightarrow 0$. Moreover, since $T$ is conformal, $ E(\mfr)  =  E(\mfr|_{A_{\nu}}) +  E(\mfr|_{B_{\nu}\backslash B_{\rho}}) = E(f_{\nu}) +  E(g_{\nu'})$ by conformal invariance of the Dirichlet integral. \medskip
     
{\it B}-II) The convergence of $\{f_{\nu}\},\,\{g_{\nu'}\}$ to $\mf^i$. \medskip

We first investigate the modulus of continuity of harmonic maps $\{h_{\nu}\}: A_{\nu} \rightarrow N$ which converge uniformly ($C^0$-norm) on $\partial B$ with $E(h_{\nu})\le L$ for some $L>0$, independent of $\nu\le \nu_0$ for some $\nu_0\in(0,1)$. We shall only deal with the assumption (C2), because the argument can clearly be applied to the case (C1) as well.\medskip

Let $G_{R}:=\overline{B_R(z)}\subset A_{\nu}$ for $\nu\le \tilde{\nu_0}$. If $z\in \partial B$, consider $G_{R}:=\overline{B_R(z)}\cap \overline{A_{\nu}}$. Given $\varepsilon>0$, by the Courant-Lebesgue Lemma there exists $\delta>0$, independent of $\nu \le \nu_0$, such that $\text{the length of}\, h_{\nu}|_{\partial G_{\delta}}$ does not exceed $\min\{\frac{\varepsilon}{4}, \frac{i(N)}{4}\}$, $i(N)>0$. Then $h_{\nu}|_{\partial G_{\delta}}\subset B(q,s)$ for some $q\in N, s\le \min\{\frac{\varepsilon}{2}, \frac{i(N)}{2}\}$. Observe that $h_{\nu}$ is continuous on $\partial G_{\delta}$, and there exists an $H^{1,2}$-extension $X$ of disc type, whose image is in $B(q,s)$ with $X|_{\partial B_{\delta}}= h_{\nu}|_{\partial B_{\delta}}$, by the same argument of Remark \ref{trace}. Thus there exists a harmonic extension $h'$ with $h'(G_{\delta}) \subset B(q,s)\subset B(q,\frac{\varepsilon}{2})$, by \cite{hkw}. From Lemma \ref{disc}, $h'$ is homotopic to $h$ on $G_{\delta}$, and from the energy minimizing property of harmonic maps, $h_{\nu}|_{G_{\delta}}=h'$. Hence, the functions $h_{\nu}$ with $\nu\le \nu_0$ have the same modulus of continuity. Furthermore, if these mappings have the same boundary image, they are $C^0$-uniformly bounded on each relatively compact domain.\medskip
   
Now apply the above result to $\{\mfr,\,\rho\le\rho_0\}$ in $\mathbb R^k$. For some $\rho_0\in(0,1)$ then, the functions $f_{\nu}$ resp. $g_{\nu'}$ have the same modulus of continuity for all $\rho \in (0,\rho_0)$, and some subsequences, denoted again by  $f_{\nu}$ resp. $g_{\nu'}$ are locally uniformly convergent.\medskip

Recall that our maps are continuous, so by localizing in both domain and image, harmonic functions, seen as solutions of Dirichlet problems, may be also regarded as weak solutions $f$ of the following elliptic systems in local coordinate charts of $N$: 
        \begin{equation}\label{ellipticpdg}
              d_{i}d_{i}f^{\alpha} = -\Gamma^{\alpha}_{\beta\gamma}d_{i}f^{\beta}d_{i}f^{\gamma}=:
                                G^{\alpha}(\cdot,f(\cdot),d f(\cdot)).
        \end{equation}
We can take the same coordinate charts for the image of $\{f_{\nu}\}_{\nu\le\nu_0}$ and $\{g_{\nu'}\}_{\nu'\le \nu'_0}$, where $\nu_0:=\nu(\rho_0), \,\nu'_0:=\nu'(\rho_0)$, to the effect that we have the same weak solution system for (\ref{ellipticpdg}). Moreover, since $h_{\alpha\beta}$ and $\Gamma^{\alpha}_{\beta\gamma}$ are smooth, the structure constants of the weak systems (see \cite{j1} section 8.5) are independent of $\rho\le \rho_0$.\medskip
      
Now consider $K_{\sigma}\dach{\sigma} = \{\sigma \le |z| \le 1-\sigma \}$, $\sigma\in (0,1)$. From the regularity theory of \cite{lu} and \cite{j1}(section 8.5) and by the covering argument, there exists $C\in \mathbb R$ such that $\left \| f_{\nu}|_{K_{\sigma}\dach{\sigma}}\right \|_{H^{4,2}} \le C \quad \text{for all}\en \nu\in (0, \nu_0)$. Hence the Sobolev's embedding theorem implies that for some sequence $\{\rho_i\} \subset (0,1)$, $\lim_{\rho_i \rightarrow 0}f_{\nu(\rho_i)}|_{K_{\sigma}\dach{\sigma}} = f' \en \text{in} \en C^2(K_{\sigma}\dach{\sigma},\mathbb R^n)$, with $\tau_h(f') = 0$ in $K_{\sigma}\dach{\sigma}$.\medskip

For $\sigma:=\frac{1}{n}$, we choose a sequence $\{f_{\nu({\rho_{n,i}})}\}$ as above such that $\{\rho_{n+1,i}\}$ is a subsequence of $\{\rho_{n,i}\}$. By diagonalizing we obtain a subsequence $\{f_{\nu({\rho_{n,n}})}\}, n\ge n_0$  which converges locally to $f'$ in the $C^2$-norm, so $f'$ is harmonic on $B\backslash(\partial B\cup \{ 0\})$.\medskip

On the other hand $f_\nu|_{\partial B} = x^1$ for all $\nu$, and the $f_\nu$'s converge uniformly to $f'$ in a compact neighbourhood of $\partial B$. Thus, $f'$ is continuous on $\overline B \backslash\{0\}$ with $f'|_{\partial B} = x^1$. Also observe that $\osc_{\partial B_r} f' \rightarrow 0$ as $r\rightarrow 0$, by construction.\medskip

For each compact $K\subset B\backslash\{ 0\}$, $\int_{K}|df'|^2d\omega= \lim_{\rho_i \rightarrow 0}\int_{K}|df_{\nu(\rho_i)}|^2 \le L$, with $L$ independent of $K$. Thus, $f'\in H\dach {1,2}(B\backslash \{0\},N)$, and $f'$ can be extended to a weakly harmonic map on $B$ (\cite{j1} Lemma 8.4.5, see also \cite{sku}, \cite{g2}). Thus, $f'$ can be considered weakly harmonic and $ f'\in C^0(\overline B,N)\cap C^2(B,N)$ with $f'|_{\partial B}=x^1$, so uniqueness forces $f' = \mf^1(x^1)$.\medskip

Similar results hold for $g_{\nu'}$.\medskip

{\it B}-III)  The convergence of the energy.\medskip

We consider  $\eta\circ f$, and denote it again by $f:=(f^{a})_{a = 1,\cdots ,k}\in H^{1,2}(\Omega,\mathbb R^k)$ for obvious reasons.\medskip

Since $\eta$ is isometric, for  $f:=(f^{\alpha})_{\alpha = 1,\cdots ,n}\in H^{1,2}(\Omega,N)$ we have $\int_{\Omega}|d(f^\alpha)|_h^2d\omega = \int_{\Omega}|d(f^a)|_{\mathbb R^k}^2d\omega$. A harmonic map $f\in H^{1,2}(\Omega,N)$ satisfies  
           \begin{equation}\label{secondfundamental}
                  \int_{\Omega}(\langle df,d\psi\rangle  - \langle II\circ f(df,df),\psi\rangle )dw = 0
           \end{equation}
for any $\psi \in H^{1,2}_0\cap C^0(\Omega,\mathbb R^k)$, where II is the second fundamental form of $\eta$.\medskip 
   
Set $K_{\sigma} = \{\sigma \le |z| \le 1\}, \sigma > 0$ and we consider $\mathbb R^k$-harmonic maps $H_\nu$  and $\widetilde{H}_\nu$ on $K_{\sigma}$ with $ H_{\nu}|_{\partial K_{\sigma}} = f_\nu|_{\partial K_{\sigma}}$ and $\widetilde{H}_{\nu}|_{\partial K_{\sigma} } = \mf^1|_{\partial K_{\sigma}}$, where $\nu \in (0,\sigma)$. Let $H : B \rightarrow \mathbb R^k$ be harmonic with $H|_{\partial B} = H_{\nu}|_{\partial B} = \widetilde{H}_{\nu}|_{\partial B} = x^1$. Then $\{H_\nu\}, \{\widetilde{H}_\nu\}$ have the same modulus of continuity up to $\partial B$, and we have $\|H_\nu - H\|_{C^0;K_{\sigma}} \rightarrow 0, \quad \|\widetilde{H}_\nu - H\|_{C^0;K_{\sigma}} \rightarrow 0 \quad \text{as} \en \nu \rightarrow 0.$ Furthermore, for $X_\nu := (f_\nu - \mf^1) + (H_\nu - \widetilde {H}_\nu)\in H^{1,2}_0\cap C^0(K_{\sigma},\mathbb R^k)$, we obtain 
        \begin{equation*}
          \|X_\nu\|_{(C^0;K_{\sigma})} \le \|f_\nu - \mf^1\|_{C^0;K_{\sigma}} +  \|H_\nu - H\|_{C^0;K_{\sigma}} + 
               \|H - \widetilde {H}_\nu\|_{C^0;K_{\sigma}}\rightarrow 0 \quad \text{as} \en \nu \rightarrow 0. 
        \end{equation*}
  
Now consider  
        \bqrs
           \lefteqn{\int_{K_{\sigma}}\langle d(f_{\nu}-\mf^1),d(f_{\nu}-\mf^1)\rangle d\omega }\hspace{1cm}\\ 
               && =  \underbrace{\int_{K_{\sigma}}\langle d(f_{\nu}-\mf^1),dX_{\nu}\rangle d\omega}_{:=I} -  
               \underbrace{\int_{K_{\sigma}}
                              \langle d(f_{\nu}-\mf^1),d(H_{\nu}-\widetilde{H}_{\nu})\rangle d\omega}_{:=II} .
        \eqrs
When $\nu \rightarrow 0$
        \bqr
           |I| 
      &\le & \left | \int_{K_{\sigma}}\langle II\circ f_{\nu}(df_{\nu},df_{\nu}),X_{\nu}\rangle d\omega \right | + 
                   \left | \int_{K_{\sigma}}\langle II\circ (d\mf^1,d\mf^1),X_{\nu}\rangle d\omega\right | \nonumber \\
               & = & C(\|f_{\nu}\|_{1,2;0}),\|\mf^1\|_{1,2;0}) \| X_{\nu}\|_{(C^0;K_{\sigma})} \rightarrow 0
              \label{CX}
        \eqr
      from (\ref{secondfundamental}). Moreover, since $H_\nu - \widetilde {H}_\nu$ is harmonic in $\mathbb R^k$,
        \bqr 
          |II| 
                \le  \int_{\partial B_{\sigma}}\left | \partial_{r}(H_\nu - \widetilde {H}_\nu)\right |d\omega 
                       \|f_{\nu} - \mf^1\|_{C^0;K_{\sigma}} \rightarrow 0 
                          \quad \text{as} \en \nu \rightarrow 0.\label{KH}
        \eqr
Thus $\int_{K_{\sigma}}|d(f_{\nu}-\mf^1)|^2 d\omega \rightarrow 0$, and $\int_{K_{\sigma}}\left | df_\nu \right | ^2d\omega \rightarrow \int_{K_{\sigma}}\left | d\mf^1\right | ^2d\omega$, for any $K_{\sigma}$. Since $ \int_{B_{\sigma}}|d\mf^1|^2d\omega \rightarrow 0$ as $\sigma \rightarrow 0$,  we obtain $\int_{A_\nu}\left | df_\nu \right | ^2d\omega \rightarrow \int_{B}\left | d\mf^1\right | ^2d\omega $ as $ \nu \rightarrow 0$. Similarly $ \int_{A_{\nu'}}\left | dg_\nu' \right | ^2d\omega \rightarrow \int_{B}\left | d\mf^2\right | ^2d\omega \quad \text{as} \en \nu' \rightarrow 0$.\medskip
     
Now to the uniform convergence on $\mathcal N_{\varepsilon}(x^i_0)$. Replace $f(\overline{\ar})$ by $\overline{B(p,r)}$(for (C1)) or $N$ (for (C2)) in the proof of Lemma \ref{tangent}. Then, $\|\mfr(x^1,x^2)\|_{H^{1,2}}\le C \en \text{uniformly on} \en \mathcal N_{\varepsilon}(x^i_0)$, where the constant $C$ depends on $x^i_0$, while $\varepsilon$ does not. The convergence in (\ref{CX}),\en (\ref{KH}) is uniform on $\mathcal N_{\varepsilon}(x^i_0)$. The proof of {\it (B)} is eventually completed.\bigskip

            
{\it C}-I)The set-up.\medskip

We must show that for $x^i\in M^i$ and $\xi^i\in \mathcal T_{x^i}$, 
        \begin{equation*}
          \langle\delta_{x^i}\me_{\rho},\xi^i\rangle  \en \longrightarrow \en \langle\delta_{x^i}\me,\xi^i\rangle\quad 
            \text{uniformly} \en \text{on}\en 
                              \mathcal N_{\varepsilon}(x^i_0)\subset M^i, i=1,2 \en \text{as} \en \rho \rightarrow 0.   
        \end{equation*}
It suffices to show the assertion for $i=1$. We know that
          \bqrs
              \lefteqn{\langle\delta_{x^1}\me_{\rho},\xi^1\rangle  
                   =  \int_{A_{\nu(\rho)}}\langle d\mfr(x^1,x^2),\nabla\mj_{\mfr}(\xi^1,0)\rangle _hd\omega 
                         + \int_{B_{\nu(\rho)}\backslash B_{\rho}}
                          \langle d\mfr(x^1,x^2),\nabla\mj_{\mfr}(\xi^1,0)\rangle _hd\omega}\hspace{2.8cm}\\ 
                 & = & \int_{A_{\nu(\rho)}}\langle d\mfr(x^1,x^2),\nabla\mj_{\mfr}(\xi^1,0)\rangle _hd\omega 
                         + \int_{A_{\nu'}}\langle dg_{\nu'},\nabla\mj_{g_{\nu'}}(0,\zeta_{\nu'})\rangle d\omega,
         \eqrs
where $g_{\nu'}(\cdot) = \mfr\circ T(\cdot) \,$ and  $\zeta_{\nu'}(\nu'e^{i\theta}) = \mj_{\mfr}(\xi^1,0)
(\nu e^{i\theta})\,$ with $\nu' :=\frac{\rho}{\nu(\rho)}$. Observe $\mj_{\mfr}(\xi^1,0)\circ T $ is a Jacobi field along $g_{\nu'}$, by the conformal property of $T$. \medskip

{\it C}-II) The convergence of Jacobi fields.\medskip

First, let $V_{\nu}:= \mj_{\mfr}(\xi^1,0)|_{A_{\nu}} = v^{\alpha}_{\nu}\frac{\partial}{\partial y^{\alpha}}\circ f_{\nu}$, for which we will show the existence of a $\nu_0\in(0,1)$ giving 
          \begin{equation}\label{dv}
                   \|DV_\nu \|_2^2 := \int _{A_\nu}h_{\alpha\beta}\circ f_\nu v_{\nu,i}^{\alpha}v_{\nu,i}^{\beta}
                      d\omega \le C \en \text{for all}\en  \nu \in(0,\nu_0).
          \end{equation}
By direct computation $\|DV_\nu \|_2^2 \le CE(V_{\nu}) + C(N,\|V_\nu\|_{C^0},\|f_\nu\|_{C^0},E(f_\nu))$. Since Lemma \ref{jacobifield} yields $\|V_{\nu}\|_{C^0} \le \|\xi^1_{\nu}\|_{C^0}$, we only need to show that
             \begin{equation}\label{dul}
                  E(V_\nu) := \int_{A_\nu}| \nabla^{f_\nu}V_\nu |^2d\omega \le C,\quad  \nu \in(0, \nu_0). 
             \end{equation}

Let $X_{\nu} := x^\alpha_{\nu}\frac{\partial}{\partial y^\alpha}\circ f_{\nu}\in H^{1,2}(A_{\nu}, f_{\nu}^{\ast}TN)$, where $ x^\alpha_{\nu}(z) := v^\alpha_{2\nu_0}(\tau^{2\nu_0}_{\nu_0}(z)), \, \nu_0 \le |z| \le 1$ (see section \ref{construction} for the definition of $\tau^{2\nu_0}_{\nu_0}$) and $x^\alpha_{\nu}(z) := 0, \,\nu \le |z| \le \nu_0$. Clearly, $\|DX_{\nu}\|^2_2 \le C(\nu_0,N) \|DV_{2\nu_0} \|_2^2$ for all $\nu\le \nu_0$. \medskip 
          
By the minimality property of Jacobi fields and Young's inequality,
               \bqrs
                 \lefteqn{\int_{A_{\nu}}(|\nabla^{f_{\nu}}(V_{\nu})|^2 
                             - \langle  trR(df_{\nu},V_{\nu})df_{\nu},V_{\nu}\rangle )d\omega 
                      \le  \int_{A_{\nu}}(|\nabla^{f_{\nu}}(X_{\nu}) |^2 
                               - \langle  trR(df_{\nu},X_{\nu})df_{\nu},X_{\nu}\rangle)d\omega}\\ 
              & \le &  \int_{A_{\nu}}\hab\circ f_{\nu}x_{\nu,i}^{\alpha}x_{\nu,i}^{\beta}d\omega
                        +\varepsilon\int_{A_{\nu}}|x_{\nu,i}^{\alpha}\frac{\partial}{\partial y^{\alpha}}\circ 
                      f_{\nu}|_h^2d\omega
                              +\varepsilon^{-1}\int_{A_{\nu}}|x_{\nu}^{\gamma}f^{\delta}_{,i}
                                  \Gamma^{\beta}_{\gamma\delta}\circ f_{\nu}
                             \frac{\partial}{\partial y^{\beta}}\circ f_{\nu}|_h^2d\omega \\
               &  & +\int_{A_{\nu}}\hab\circ f_{\nu}x_{\nu}^{\gamma}x_{\nu}^{\lambda}f^{\delta}_{,i}f^{\mu}_{,i}
                      \Gamma^{\alpha}_{\gamma\delta}\circ f_{\nu}\Gamma^{\beta}_{\lambda\mu}\circ f_{\nu}d\omega 
                              - \int_{A_{\nu}}\langle trR(df_{\nu},X_{\nu})df_{\nu},X_{\nu}\rangle d\omega\\
              & \le & C(N,\varepsilon,\| f_{\nu}\|_{C^0},E(f_{\nu}),\|V_{2\nu_0}\|_{C^0},\|DV_{2\nu_0}\|_2^2).
             \eqrs
But $E(V_\nu)\le C,\en \nu \in(0, \nu_0)$, since 
           \[\int_{A_{\nu}}\langle trR(df_\nu,V_\nu)df_\nu,V_\nu\rangle )d\omega       
                                 \le C(N,\| f_{\nu}\|_{C^0},E(f_{\nu}),\|\xi^1\|_{C^0}).\]
Therefore we have (\ref{dv}), and this means that $\{(v^{\alpha}_\nu)|\nu\le\nu_0 \}_{\alpha = 1,\cdots ,n}$ has the same modulus of continuity, see the argument in {\it B}-III) and Lemma \ref{jacobifield}.\medskip

With the same charts as in {\it (B)}, $(v\dach{\alpha}_{\nu(\rho)})\in \mathbb R\dach n,\en \nu\le \nu_0$ are weak solutions of the Jacobi fields system with uniformly bounded energy and same modulus of continuity on $K_{\sigma} = \{\sigma \le |z| \le 1\}$, with $\sigma >0$ for small $\rho$, again by Lemma \ref{jacobifield}. Just as in {\it (B)}, $\{V_{\nu}\}$ converges to the Jacobi field along $\mf\dach 1|_{B\backslash \{0\}}$ with boundary $\xi^1$, and for $\mj_{\mf^1}(\xi^1)=:w^{\alpha}\frac{\partial}{\partial y^{\beta}}\circ \mf^1$, we have 
       \[\|(v^{\alpha}_\nu(z))-(w^{\alpha}(z))\|_{C^0;K_{\sigma}}\rightarrow 0, \en 
          \|(v^{\alpha}_\nu(z))-(w^{\alpha}(z))\|_{C^2;K}\rightarrow 0, \quad \text{as}\en  \nu(\text{or}\en \rho)
          \rightarrow 0,\]
on any compact $K\subset B\backslash \{0\}$.\medskip

{\it C}-III) The convergence of derivatives.\medskip

Taking $K_{\sigma}$ as above, we denote $f_{\nu}|_{K_{\sigma}}$ and $\mf^1|_{K_{\sigma}}$ by $f_{\nu}$ and $\mf^1$, respectively.\medskip
 
Note that $\exp_{\mf^1}:\mathcal U(0) \rightarrow H^{1,2}\cap C^0(K_{\sigma},N)$ is a diffeomorphism on some neighbourhood $\mathcal U(0)\in H^{1,2}\cap C^0(K_{\sigma},(\mf^1)^{\ast}TN)$, because $d(\exp_{\mf^1})_0 = Id$. Moreover, $\| f_{\nu}- \mf^1|_{K_{\sigma}}\|_{H^{1,2}\cap C^0} \rightarrow 0$ as $\nu \rightarrow 0,$ so there exists $\xi_{\nu}\in H^{1,2}\cap C^0(K_{\sigma},(\mf^1)^{\ast}TN)$ for small $ \nu>0$ with $\exp_{\mf^1}\xi_{\nu} = f_{\nu}$.\medskip
 
The mapping $\xi\mapsto d\exp_{\mf^1,\xi}$ depends smoothly on $\xi_\nu\in T_{\mf^1}H^{1,2}\cap C^0(K_{\sigma},N)$, so $d\exp_{\mf^1,\xi_{\nu}} \rightarrow Id \en \text{in} \en H^{1,2}\cap C^0(K_\sigma)$, since $\xi_\nu \rightarrow 0$ in $H^{1,2}\cap C^0(K_{\sigma},(\mf^1)^{\ast}TN) \en \text{as}\en  \nu \rightarrow 0$. For $ W_\nu := w^{\alpha}_{\nu}\frac{\partial}{\partial y^{\alpha}}\circ \mf^1 := d\exp_{\mf^1,\xi_{\nu}}^{-1}(V_\nu)$ we have $ \|w^{\alpha}_\nu(z)-w^{\alpha}(z)\|_{C^0;K_{\sigma}}\rightarrow 0$ by {\it C}-II). Moreover, $d\mf^1 \rightarrow df_{\nu}$ in $L^2$, thus $\int_{K_{\sigma}}|d\exp_{\mf^1,\xi_{\nu}}(d\mf^1) - df_{\nu}|^2 d\omega \rightarrow 0$.\medskip
          
We next observe, for $\nabla^{\mf^1}W_{\nu}=(w_{\nu,i}^{\alpha} + w_{\nu}^{\gamma}(\mf^1)^{\beta}_{,i}\Gamma^{\alpha}_{\beta\gamma}(\mf^1))dz^i\otimes \frac{\partial}{\partial y^{\alpha}}\circ\mf^1$, that   
           \bqr\label{thirty}
              \int_{K_{\sigma}}|d\exp_{\mf^1,\xi_{\nu}}(\nabla^{\mf^1}W_{\nu})- 
                                    \nabla^{f_{\nu}}V_{\nu}|^2 d\omega 
                           \rightarrow 0 \en \text{as} \en \nu\rightarrow 0,
           \eqr
since $\|\mf^1 - f_{\nu}\|_{1,2;0}\rightarrow 0$, $d\exp_{\mf^1,\xi_{\nu}} \rightarrow Id$ in $C^0$, $\partial_i(d\exp_{\mf^1,\xi_{\nu}}) \rightarrow \partial_i (Id) = 0$ in $L^2$.\medskip    
               
Thus, for $X_{\nu}, Y_{\nu} \in H^{1,2}\cap C^0(K_{\sigma},T^{\ast}M\otimes f_{\nu}^{\ast}TN)$ with $\int_{K_{\sigma}}|X_{\nu}|^2d\omega \rightarrow 0, \en  \int_{K_{\sigma}}|Y_{\nu}|^2d\omega  \rightarrow 0$, 
             \[d\exp_{\mf^1,\xi_{\nu}}(d\mf^1) = df_{\nu} + X_{\nu}, 
                 \en d\exp_{\mf^1,\xi_{\nu}}(\nabla^{\mf^1}W_{\nu}) =  \nabla^{f_{\nu}}V_{\nu}+ Y_{\nu}.\] 
Gau{\ss} lemma prescribes that $ \langle d\mf^1, \nabla^{\mf^1}W_{\nu} \rangle_h = \langle df_{\nu} + X_{\nu}, \nabla^{f_{\nu}}V_{\nu} + Y_{\nu}\rangle_h$. Thus the H\"older inequality and (\ref{dul}) give
          \bqr
            \lefteqn{ \int_{K_{\sigma}} \left(\langle df_{\nu},\nabla ^{df_{\nu}}V_{\nu}\rangle_h - 
                       \langle d\mf^1,\nabla^{\mf^1}\mj_{\mf^1}(\xi^1)\rangle_h \right)d\omega}\hspace{2.0cm}\nonumber\\
                & = & \int_{K_{\sigma}} \left(\langle d\mf^1,\nabla ^{\mf^1}W_{\nu}\rangle_h - 
                          \langle d\mf^1,\nabla^{\mf^1}\mj_{\mf^1}(\xi^1)\rangle_h \right)d\omega + o(1)\nonumber\\
               & \le & E(d\mf^1)\| \nabla ^{\mf^1}W_{\nu} - \nabla^{\mf^1}\mj_{\mf^1}(\xi^1)\|_{L^2;K_{\sigma}} + o(1).
                                 \label{edf}
           \eqr
In order to estimate the last term, consider $A_\nu:=a^{\alpha}_\nu\frac{\partial}{\partial y^{\alpha}}\circ \mf^1$ and $A:=a^{\alpha}\frac{\partial}{\partial y^{\alpha}}\circ \mf^1$ such that $d\eta\big(a^{\alpha}_\nu\frac{\partial}{\partial y^{\alpha}}\circ \mf^1\big)$ and $d\eta\big(a^{\alpha}\frac{\partial}{\partial y^{\alpha}}\circ \mf^1\big)$ are harmonic in $\mathbb R^k$ with $A_\nu|_{\partial K_{\sigma}}=W_\nu|_{\partial K_{\sigma}}$, $A|_{\partial K_{\sigma}}=W|_{\partial K_{\sigma}}$, for $W:= \mj_{\mf^1}(\xi^1)$. Clearly, $\|d\eta(A_\nu-A)\|_{1,2;0} \rightarrow 0$.\medskip

Now, consider a test vector field $Z_{\nu}:=W_\nu- W - A_\nu + A \in H^{1,2}_0\cap C^0(K_{\sigma},(\mf^1)^{\ast}TN)$. Observing that  $W$ and $V_\nu$ are Jacobi fields along $\mf^1|_{K_\sigma}$ and $f_\nu|_{K_\sigma}$ respectively, we have 
         \bqrs
              \lefteqn{\int_{K_\sigma}\langle \nabla^{\mf^1}(W_{\nu} - W),
                   \nabla^{\mf^1}Z_{\nu}\rangle_hd\omega}\label{estimate}\hspace{0.4cm}\\
              & = &\int_{K_\sigma}\{ \langle \nabla^{\mf^1}W_{\nu},
                      \nabla^{\mf^1}Z_{\nu}\rangle_h  - \langle trR\circ \mf^1(W,d\mf^1)d\mf^1,
                         Z_{\nu}  \rangle_h \nonumber\\
               & & \hspace{0.9cm}  -\langle \nabla^{f_\nu}V_\nu,
                         \nabla^{f_\nu}(L_\nu (Z_{\nu})) \rangle_h +\langle trR\circ f_\nu(V_\nu,df_\nu)df_\nu,
                         (L_\nu(Z_{\nu}))\rangle_h \} d\omega\nonumber\\
              & = &\int_{K_\sigma}\{ \langle \nabla^{\mf^1}W_{\nu},
                          \nabla^{\mf^1}Z_{\nu}\rangle_h  - \langle trR\circ \mf^1(W,d\mf^1)d\mf^1,
                                               Z_{\nu}\rangle_h\nonumber \\ 
               & &  \hspace{0.2cm} -\langle \nabla^{\mf^1}L^{-1}_\nu(V_\nu),
                               \nabla^{{\mf^1}}Z_\nu  \rangle_h  +\langle trR\circ f_\nu(V_\nu,df_\nu)df_\nu,
                                  (L_\nu (Z_{\nu}))\rangle_h \} d\omega + o(1) \nonumber
         \eqrs
with $L_\nu:=d\exp_{\mf^1,\xi_{\nu}}$. This expression converges to $0$ as $\nu\rightarrow 0$, since $L_\nu^{-1}(V_\nu)= W_\nu$, $ \| Z_{\nu} \|_{C^0;K_{\sigma}}\rightarrow 0$ and $ \|\mf^1\|_{1,2;0}, \|W\|_{C^0}, \|f_\nu\|_{1,2;0}, \|V_\nu\|_{C^0} < C \en \text{for all}\en \nu \in(0,\nu_0)$.\medskip
   
Moreover, $\int_{K_\sigma}|\nabla^{\mf^1}(A_{\nu} - A)|_h^2 d\omega \rightarrow 0$, since $\|d\eta(A_\nu-A)\|_{C^0}\rightarrow 0$ and because of (\ref{yeoseot}). Thus, (\ref{edf}) converges to $0$ for each $\sigma \in (0,1)$. Now let $\sigma \rightarrow 0$. Then  
        \[\int_{A_{\nu(\rho)}}\langle d\mfr(x^1,x^2),\nabla\mj_{\mfr}(\xi^1,0)\rangle_hd\omega \rightarrow 
                \int_{B}\langle d\mf^1(x^1),\nabla\mj_{\mf^1}(\xi^1)\rangle_hd\omega, \quad \rho \rightarrow 0,\]
since $\int_{B_{\sigma}}\langle d\mf^1(x^1),\nabla\mj_{\mf^1}(\xi^1)\rangle_hd\omega \rightarrow 0$ as $\sigma \rightarrow 0$.\medskip

In a similar way, $\int_{A_{\nu'(\rho)}}\langle dg_{\nu'},\nabla\mj_{g_{\nu'}}(0,\zeta_{\nu'})\rangle d\omega \rightarrow\int_{B}\langle d\mf^2(x^2),\nabla\mj_{\mf^2}(0)\rangle_hd\omega = 0$.\medskip

The uniform convergence on $\mathcal N_{\varepsilon}(x^i_0)$ is clear.\medskip
       
In this manner we could also show that $\delta_{x^1}\me_{\rho}, \delta_{x^2}\me_{\rho}$ are continuous with respect to $\rho \in (0,1)$, and uniformly continuous on  $\mathcal N_{\varepsilon}(x^i_0)$. This concludes part {\it (C)}.\bigskip

{\it (D)} Along the lines of \cite{s4}, the differential form 
         \bqr
          \frac{\partial}{\partial t}|_{t=\rho} \me(x^1,x^2,t) 
           =\int^{2\pi}_0\int^{1}_{\rho} \left[ |\partial_r \mfr|^2 - \frac{1}{r^2}|\partial_{\theta}\mfr|^2 \right]
                  \frac{1}{1-\rho} drd\theta \label{g3} 
            \eqr                 
proves {\it (D)}, bringing to an end the proof of Lemma \ref{differ}.  \hfill $\Box$


\subsection{Critical points of $\me$}\label{threehundred}

For given Jordan curves $\Gamma_1, \Gamma_2, \Gamma$ in $(N,h)$ with $\text{dist}(\Gamma_1,\Gamma_2)>0$, we consider the Plateau problems $\mathcal P(\Gamma_1,\Gamma_2)$ and $\mathcal P(\Gamma)$.\medskip
 
We define for $x=(x^1,x^2,\rho)\in \overline{\mm}$, 
    \bqr
          g_i(x) &:= &\sup\limits_{\begin{array}{c}\xi^i\in \mathcal T_{x^i} \\ \|\xi^i\| < l_i \end{array}}
                                      (-\langle \delta_{x^i}\me,\xi^i\rangle),\quad\quad i=1,2,\label{critical1} \\
             g_3(x) & := & \left\{ \begin{array}{c@{\quad , \quad}l}
                              \left|\rho\cdot \partial_{\rho}\me\right| & \rho > 0 \\
                                                                                  0 & \rho = 0, 
                         \end{array} \right.\nonumber\\
           g(x) & := & \Sigma _{j=1}^{3} g_j(x).\nonumber
    \eqr
In the definition of $l_i$ of section \ref{construction}, we can clearly require that $l_i\le \{1,i_{\tilde{h}}(\Gamma_i)\}$. Note that $g_j\ge 0, j=1,2,3$, because $g_i(x)<0,\en i=1,2$ would imply $\langle \delta_{x^i}\me,\xi^i\rangle \ge \sigma >0$ for all $\xi^i \in \mathcal T_{x^i}$ with $\|\xi^i\| < l_i$. Since $\mathcal T_{x^i}$ is convex, $\langle \delta_{x^i}\me,t\xi^i\rangle = t\sigma \ge\sigma$, $t\in [0,1]$, a contradiction. Clearly, $g_3(x)\ge 0$. Now we are ready to define the critical points of $\me$.\medskip

{\bf Definition}  $x\in \overline{\mm}$ is a critical point of $\me$ if $g(x) =0$, i.e. $g_j=0, j=1,2,3$.

\begin{lemma}\label{gi}
The functions $g_j$ are continuous, $j=1,2,3$. In particular, as $\rho \rightarrow \rho_0\in[0,1)$, $g_j(x^1,x^2,\rho)$ converges uniformly to $g_j(x^1,x^2,\rho_0)$ on $\mathcal N_{\varepsilon}(x^i)$, $i=1,2$, for some small $\varepsilon > 0$.
\end{lemma}
{\bf Proof.} The uniform convergence of $g_i$ follows immediately from the uniform convergence of $\delta_{x^i}\me$, see Lemma \ref{differ} {\it(C)}.\medskip
    
Let $\{x_n\}=\{(x^1_n, x^2_n, \rho_n)\} \subset \overline{\mm}$ strongly converge to $x=(x^1, x^2, \rho)$. From the above, $g_i(x^1_n, x^2_n, \rho_n)\rightarrow g_i(x^1_n, x^2_n, \rho)\en \text{uniformly on} \en \{n\ge n_0\}$.\medskip
   
Let $\tilde{x}_n:=(x^1_n, x^2_n, \rho)$ and $\widetilde{\exp} _{x^i_n}\xi^i_n = x^i$. Observe that $d\widetilde{\exp}_{x^i_n,\xi^i_n}\rightarrow Id $ in $\hs\cap C^0$, hence  for some $t_0$ independent of  $n\ge n_0$, $\|t_0d\widetilde{\exp}_{x^i_n,\xi^i_n}(\eta^i_n)\|_{\mathcal T_{x^i}}<l_i \en\text{if} \en \|\eta^i_n\|_{\mathcal T_{x^i_n}} < l_i$. Note that $\mathcal T_{x^i}$ is convex and contains zero.\medskip

Then by Lemma \ref{differ} {\it(A)}, for given $\delta>0$ there exist $t_0(\delta)$ and $n_0(\delta)$ as above such that for each  $\|\eta^i_n\|_{\mathcal T_{x^i_n}} < l_i$ with $n\ge n_0(\delta)$,
     \bqrs  
      -\langle \delta_{x^i}\me(\tilde{x}_n), \eta^i_n\rangle & \le & -\langle \delta_{x^i}\me(x), 
                                                   d\widetilde{\exp}_{x^i_n,\xi^i_n}(\eta^i_n)\rangle + \delta \\
        & \le & -\langle \delta_{x^i}\me(x), t_0d\widetilde{\exp}_{x^i_n,\xi^i_n}(\eta^i_n)\rangle + 2\delta
         \le  g_i(x) + 2\delta.
   \eqrs
This implies $g_i(\tilde{x}_n) \le g_i(x) + 2\delta$. On the other hand $g_i(x) \le g_i(\tilde{x}_n) + 2\delta$, so $g_i(x^1_n, x^2_n, \rho) \rightarrow g_i(x^1, x^2, \rho)$ as $n\rightarrow \infty$.\medskip 
   
Together with the above uniform convergence on $\mathcal N_{\varepsilon}(x^i)$ for $\rho_n\rightarrow \rho$, we infer the continuity of $g_i$, $i=1,2$. The continuity and uniform continuity of $g_3$ are easy consequences of the expression of $\frac{\partial}{\partial \rho}\me$. \hfill $\Box$ 

\begin{proposition}\label{h22}
   $x=(x^1,x^2,\rho)\in M^1\times M^2\times [0,1)$ is a critical point of $\me$ if and only if 
 $\mf_\rho(x^1,x^2)$ (for $\rho \in  (0,1))$, resp. $\mf^i(x^i)$ is a solution of $\mathcal P(\Gamma_1,\Gamma_2)$,  
 resp. $\mathcal P(\Gamma_i), i=1,2$.
\end{proposition}
{\bf Proof.}   (I) Let $x=(x^1,x^2,\rho)\in M^1\times M^2\times [0,1)$ be a critical point of $\me$. From \cite{hkw} $\mf$ is continuous up to the boundary. We must show that $\mfr(x^1,x^2)$(for $\rho>0$) and $\mf^i(x^i)$ are conformal. We will show this only for $\mfr(x^1,x^2)$, the other case being analogous.\medskip

For $x\in \mm$ critical point of $\me$, we have $\mf_\rho(x^1,x^2)\in H^{2,2}(A_\rho , \mathbb R^k)$ from Theorem \ref{appendixthm}. The condition $\gamma^i\in C^3, i=1,2$ will be essential. Taking $\xi^1\in \mathcal T_{x^1}$, and denoting $\mf_\rho(x^1,x^2)$ and $\mj_{\mf_{\rho}}(\xi^1,0)$ by $\mfr$ and $\mj_{\rho}$, we compute $\xi^1\in \mathcal T_{x^1}$,  
     \bqr
        \langle \delta_{x^1}\me,\xi^1\rangle 
        & = &\int_{A_\rho} \langle d\mf_{\rho},
             \nabla_{\frac{d}{dt}}d\mf_{\rho}(\widetilde{\exp}_{x^1}t\xi^1,x^2)\big|_{t=0}\rangle_h d\omega\nonumber
                   = \int_{A_\rho} \langle \frac{\partial}{\partial z^i}\mf_{\rho},
                     \nabla_{\frac{\partial}{\partial z^i}}\mj_{\mf_{\rho}}(\xi^1,0)\rangle_h d\omega\nonumber\\
        & = & \int_{A_\rho}\text{div}(\langle \frac{\partial}{\partial z^1}\mf_{\rho},\mj_{\mf_{\rho}}(\xi^1,0)\rangle_h,
                   \langle \frac{\partial}{\partial z^2}\mf_{\rho},\mj_{\mf_{\rho}}(\xi^1,0)\rangle_h)d\omega
                \en(\text{since}\en\nabla_{\frac{\partial}{\partial z^i}} \frac{\partial}{\partial z^i}\mf_{\rho}=0)
                                                                                                      \nonumber\\   
        & = & \int_{\partial B}\langle \frac{\partial}{\partial z^1}\mf_{\rho}\vec n ,\xi^1\rangle_h d\omega.
                         \label{div3}
     \eqr
The work in \cite{s1} leads to the conformal property of $\mfr$.\medskip

(II) Let $\mf:= \mfr(x)$(resp. $\mf^i(x^i)$) be a minimal surface of annulus (resp. disc) type. By \cite{hh}, $\mf\in C^1(\overline{A_{\rho}},N)$ (resp. $C^1(\overline{B},N))$. Conformality implies $\frac{d\mf}{d\vec n}\cdot \frac{d}{d\theta} x^i \equiv 0$, and (\ref{div3}) says that $g_1(x)=0, g_2(x)=0$. That $g_3(x)=0$ follows from using (\ref{g3}) as well.                                     \hfill $\Box$


\section{Unstable minimal surfaces}\label{CPT}


\subsection{The Palais-Smale condition}\label{Palais}

By the conformal invariance of the energy function $E$, the Palais-Smale (PS) condition cannot be satisfied in the former setting for $\me$ (cf. \cite{s1} Lemma I.4.1). Hence we need the normalization used in \cite{s4}: With $P^i_k\in\Gamma_i$ fixed, $k=1,2,3$, $i= 1,2$, let 
      \[M^{i\ast}=\{x^i\in M^i: x^i(\cos \frac{2\pi (k-1)}{3}, \sin \frac{2\pi (k-1)}{3})=P^i_k\in\Gamma_i, 
                      \en k = 1,\,2,\,3\}.\]
Now define
     \bqrs
               \mm^{\ast} & = & \{x=(x^1,x^2,\rho)\in\mm : x^1(1,0) = P^1_1\in \Gamma_1\}, \\
       \partial \mm^{\ast}  & = & \{x=(x^1,x^2,0)\in \partial\mm : x^i\in M^{i\ast} \}.
     \eqrs
Given $x\in \mm^{\ast}$ and $x \in \partial \mm^{\ast}$ we take the variations from $\mathcal T_x\mm=\mathcal T_{x^1}\times \mathcal T_{x^2}\times \mathbb R$ and $\mathcal T_x\partial\mm = \mathcal T_{x^1}\times \mathcal T_{x^2}$ respectively, namely we use the original tangent spaces.\medskip

We consider the following topology:
\begin{itemize}
\item 
   A neighbourhood $\mathcal U_{\varepsilon}(x_0)$ of $x_0=(x^1_0,x^2_0,0) \in \partial \mm^{\ast}$ consists of all $x=(x^1,x^2,\rho)\in \overline{\mm^{\ast}}$ such that $\rho < \varepsilon$ and for each $i=1,2$, $ \inf_{\{\text{all}\en \sigma\}}\|\mf^i(x^i)\circ \sigma - \mf^i(x^i)\|_{1,2} < \varepsilon $, where $\sigma$ is a conformal diffeomorphism of $B$.
\item
   A sequence $\{x_n = ( x^1_n, x^2_n,\rho_n)\} \subset \overline\mm^{\ast}$ converges strongly  to $x = ( x^1, x^2,0)\in \partial\mm^{\ast}$, if all but finitely many $x_n$ lie in $\mathcal U_{\varepsilon}(x)$, for any $\varepsilon>0$.
\end{itemize}

With respect to this topology $g_j, j = 1,2,3,$ are continuous and uniformly continuous as $\rho\rightarrow \rho_0\in [0,1)$ on some $\varepsilon$-neighborhood of $(x^1,x^2)$, because of Lemma \ref{gi} and the invariance of the Dirichlet  integral under conformal changes.

\begin{proposition}[Palais-Smale condition]\label{psc}
Suppose $\{x_n\}$ is a sequence in $\overline{\mm^{\ast}}$ such that $\me(x_n) \rightarrow  \beta,\en g(x_n)\rightarrow 0$, as $n\rightarrow \infty$. Then there exists a subsequence of $\{x_n\}$ which converges strongly to a critical point of 
 $\me$ in $ \overline{\mm^{\ast}}$. 
\end{proposition}
{\bf Proof.} We prove this for the case $\{x_n\} \subset \mm^{\ast}$ with $0 < \rho_n <1$, $\me(x_n) \rightarrow \beta$, $g_j(x_n)\rightarrow 0$. If $\{x_n\} \subset \partial \mm^{\ast}$, the proof is similar. We may additionally suppose that $\rho_n \rightarrow \rho$.\medskip

Note that $\rho$ cannot be $1$, i.e. $0\le\rho < 1$, because for any $x=(x^1,x^2,\rho)\in \mathcal M$, $\frac{\rho}{1-\rho} \le c\me(x)$, since $ 0 <\dist(\Gamma_1, \Gamma_2)$. More on this can be found in \cite{s4}, Lemma 4.10.\medskip

Clearly $\int_{A_{\rho}}|d\eta\circ \mfr(x^1,x^2)|^2d\omega \ge \int_{A_{\rho}}|d\mhr(x^1,x^2)|^2d\omega \ge C(\rho) \Sigma_{i}\int_{B}|d\mh(x^i)|^2d\omega$. Thus Proposition II.2.2 of \cite{s1} guarantees that for some subsequence $\{w^i_n\}$ with $\gamma^i(w^i_n)=x^i_n$, we either have $\|w^i_n -  w^i\|_{C^0} \rightarrow 0$ with $\gamma^i\circ w^i\in \hs\cap C^0(\partial B,\Gamma_i)$, or $x^i_n=\gamma^i\circ w^i_n\rightarrow\text{const.}=a_i\in \Gamma_i$ in $L^1(\partial B)$. Therefore we have to distinguish four main cases, each divided in sub-steps.\medskip

{\bf (case 1)} Let $\rho\in(0,1)$ and  $\|w^i_n - w^i\|_{C^0} \rightarrow 0$, i.e. $\|x^i_n - x^i\|_{C^0} \rightarrow 0, x^i\in \hs\cap C^0,\en i=1,2$.\medskip

First,  $\gamma^i( w^i_n(\theta)) - \gamma^i( w^i(\theta)) = \underbrace{d\gamma^i( w^i_n(\theta))( w^i_n(\theta)- w^i(\theta))}_{=:I^i_n} - \underbrace{\int^{ w^i_n(\theta)}_{ w^i(\theta)}\int^{ w^i(\theta)}_{ w'}d^2\gamma^i( w'')d w''d w'.}_{=:II^i_n}$ In addition, $ \int_{A_{\rho}}|d\mhr(II^1_n,II^2_n)|^2d\omega \le C(\rho)(\|\mh(II^1_n)\|_{1,2;0}^2 +  \|\mh(II^2_n)\|_{1,2;0}^2) \rightarrow 0$, as $n \rightarrow \infty$, since $\|II^i_n\|_{\frac{1}{2},2;0} \le C\| w^i_n -  w^i\|_{C^0} (| w^i_n|_{\frac{1}{2}} + | w^i|_{\frac{1}{2}})$ by \cite{s3} (3.9).\medskip

Let $\mh_n:=\mhr(x^1_n,x^2_n), \en \mh:=\mhr(x^1,x^2),\en \mf_n:=\mfr(x^1_n,x^2_n):\ar\rightarrow N(\hookrightarrow \mathbb R^k)$.\medskip 

Since $\mh_n - \mh$ is harmonic on $\mathbb R^k$ and $\int_{A_{\rho}}\langle d\mh,d(\mh_n - \mh)\rangle d\omega = o(1)$ as $ n\rightarrow \infty$,       
     \bqrs 
        \int_{\ar}|d(\mh_n-\mh)|^2d\omega 
                =  \int_{\ar}\langle d\mf_n, d(\mh_n - \mh)\rangle d\omega + o(1)
                =  \int_{\ar}\langle d\mf_n, d(\mhr(I^1_n, I^2_n))\rangle d\omega + o(1).
      \eqrs 
Now consider $\xi^i_n :=-I^i_n \in \mathcal T_{x^i_n}$, and set $\mj^1_n:=\mj_{\mfr}(\xi^1_n,0),\en\mj^2_n:=\mj_{\mfr}(0,\xi^2_n)$. Then 
    \bqrs   
      \lefteqn{\int_{\ar}\langle d\mf_n,d\mhr(I^1_n,I^2_n)\rangle d\omega
         = \int_{\ar}\langle d\mf_n, d\mhr(I^1_n,0)\rangle d\omega + \langle d\mf_n,d\mhr(0,I^2_n)\rangle d\omega}\\
         &=& \int_{\ar}-\langle d\mf_n,d\mj^1_n\rangle d\omega 
              + \int_{\ar}\langle II\circ\mf_n(d\mf_n,d\mf_n),\mhr(I^1_n,0)+\mj^1_n\rangle d\omega \\
         & &\hspace{1cm} + \int_{\ar}-\langle d\mf_n,d\mj^2_n\rangle d\omega
              + \int_{\ar}\langle II\circ\mf_n(d\mf_n,d\mf_n), \mhr(0,I^2_n)+\mj^2_n \rangle d\omega \\
        &\le& g_i(x^1_n,x^2_n,\rho)\|\xi^i_n\|_{\frac{1}{2},2;0} +C(\|\mf_n\|_{1,2;0})\|\xi^i_n\|_{C^0}\\
        &\le& Cg_i(x_n)\|\xi^i_n\|_{\frac{1}{2},2;0} +C(\|\mf_n\|_{1,2;0})\|x^i_n-x^i\|_{C^0},
     \eqrs
where $C$ is independent of $n \ge n_0$, for some $n_0$. This follows from the observation (Remark \ref{trace}, Remark \ref{hundredtwo} and Lemma \ref{gi}) that $\|x^i_n - x^i\|_{C^0} \rightarrow 0$ implies the uniform convergence of $g_i(x^1_n,x^2_n,\rho_{n'}))$ on $\{x^i_n|n\ge n_0\}$ as $\rho_{n'}\rightarrow \rho$. Moreover $\|\xi^i_n\|$ are uniformly bounded.\medskip
 
Therefore $\int_{\ar}|d(\mh_n-\mh)|^2d\omega \rightarrow 0$, and $x^i_n \rightarrow x^i \en \text{strongly in} \en \hs\cap C^0(\partial B,\mathbb R^k)$.\medskip
   
{\bf (case 2)} Let $\rho\in(0,1)$,  $ \|x^1_n - x^1\|_{C^0} \rightarrow 0, \en  x^2_n =\gamma_2\circ w^2_n\rightarrow\text{const.}=a_2\in \Gamma_2 \en \text{in} \en L^1(\partial B,\mathbb R^k)$.\medskip

I) We first claim that $\mf:=\mfr(\gamma^1\circ w^1, a^2)$ is well defined and conformal. The proof is split into four steps I-a) {\bf ---} I-d) . \medskip
    
I-a) Let $x^2_n := \gamma_2\circ w^2_n$, $ a_2 := \gamma^2\circ w^2$ and $\mf_{\rho_n}:=\mf_{\rho_n}(x^1_n, x^2_n)$.\medskip
      
There must exist $\theta_0 \in [0,2\pi](\cong \partial B)$ such that $\big|\lim_{\theta\rightarrow \theta_0+} w^2(\theta) - \lim_{\theta\rightarrow \theta_0-}w^2(\theta)\big|=2\pi$. By the Courant-Lebesgue Lemma,  for given $\varepsilon >0$ there exists $r_n\in (\delta, \sqrt{\delta})$ for small $ \delta:=\delta(\varepsilon)>0$ such that with $B_{r_n}:= B_{r_n}(\theta_0)\subset\mathbb R^2$ 
     \begin{equation}\label{osc}
        \text{osc}_{A_{\rho_n}\cap \partial B_{r_n}} \mf_{\rho_n}(x^1_n,x^2_n)
         \le C\,\frac{\me(x^1_n,x^2_n,\rho_n)}{\ln(\delta^{-1})}  \le \frac{C}{\ln(\delta^{-1})} <\varepsilon.
        \end{equation}
For $\varepsilon:= \frac {1}{n}$, $C^2_n:=\partial B_{\rho_n}\backslash B_{r_n}\cup (A_{\rho_n}\cap \partial B_{r_n})$, $ Y^2_n := \mf_{\rho_n}(C^2_n)$ we see that $\text{dist}(Y^2_n, a_2) \rightarrow 0 \en \text{as} \en n\rightarrow \infty$, and the energy of $\mf_{\rho_n}|_{C^2_n}$ converges to $0$.\medskip
      
I-b) Let $\mh_{\rho_n}:=\mh_{\rho_n}(x^1_n, x^2_n)$, $\tilde{\mh}_n:=\mh_{\rho_n}(x^1, a^2)$, $\mf_{\rho_n}:= \mf_{\rho_n}(x^1_n, x^2_n)$. As above, we can say 
     \bqrs 
      \int_{A_{\rho_n}\backslash B_{r_n}}
               |d(\mh_{\rho_n}-\tilde{\mh}_n)|^2d\omega 
       & = & \int_{A_{\rho_n}\backslash B_{r_n}}
              \langle d\mf_{\rho_n}, d(\mh_{\rho_n}-\tilde{\mh}_n)\rangle d\omega + o(1)\\
       & = & \int_{A_{\rho_n}\backslash B_{r_n}}
              \langle d\mf_{\rho_n}, 
                 d\mathcal K_{\rho_n}(I^1_n, \mf_{\rho_n}|_{C^2_n}-a^2)\rangle 
                            d\omega + o(1),
    \eqrs 
where $K_{\rho_n}(I^1_n, \mf_{\rho_n}|_{C^2_n}-a^2) : A_{\rho_n}\backslash B_{r_n}\rightarrow \mathbb R^k$  denotes the Euclidean harmonic extension with $I^1_n$ on $\partial B$ and $\mf_{\rho_n}|_{C^2_n}-a^2$ on $C^2_n$.\medskip 

Let $\tilde {\mj}_n  :=\mj_{\mf_{\rho_n}}(\xi^1_n,0)$ with $\xi^1_n := -I^1_n$ and $l_n := \tilde {\mj}_n|_{C^2_n}$. Since $\|x^1_n - x^1\|_{C^0} \rightarrow 0$, it follows $\|I^1_n\|_{C^0}\rightarrow 0$ as $n\rightarrow \infty$. We can then estimate further  
        \bqrs
          \lefteqn{\int_{A_{\rho_n}\backslash B_{r_n}}
                       \langle d\mf_{\rho_n}, 
                 d\mathcal K_{\rho_n}(I^1_n, \mf_{\rho_n}|_{C^2_n}-a^2)\rangle 
                            d\omega \hspace{1.5cm}}\\
           & = & \int_{A_{\rho_n}\backslash B_{r_n}}
                    -\langle d\mf_{\rho_n}, d\tilde{\mj}_n \rangle d\omega
               + \int_{A_{\rho_n}\backslash B_{r_n}}
                     \langle d\mf_{\rho_n}, d\mathcal K_{\rho_n}(I^1_n,-l_n)
                                              +d\tilde {\mj}_n \rangle d\omega\\ 
            & & \hspace{3.5cm}   + \int_{A_{\rho_n}\backslash B_{r_n}}
                     \langle d\mf_{\rho_n}, d\mathcal K_{\rho_n}
                                                    (0, l_n + \mf_{\rho_n}|_{C^2_n}-a^2)
                                              \rangle d\omega\\
           & = & \int_{A_{\rho_n}\backslash B_{r_n}}
                    -\langle d\mf_{\rho_n}, d\tilde {\mj}_n \rangle d\omega
               + \int_{A_{\rho_n}\backslash B_{r_n}}
                    \langle II\circ \mf_{\rho_n} (d\mf_{\rho_n},d\mf_{\rho_n}), 
                                        \mathcal K_{\rho_n}(I^1_n, -l_n)+ \mj^1_n 
                                               \rangle d\omega  + o(1)\\
            & & (\text{observing that}\en 
                    \int_{A_{\rho_n}\backslash B_{r_n}}
                        \langle d\mf_{\rho_n}, d\mathcal K_{\rho_n}
                          (0, l_n + \mf_{\rho_n}|_{C^2_n}-a^2)\rangle d\omega=o(1))\\ 
           & = & \int_{A_{\rho_n}\backslash B_{r_n}}
                   -\langle d\mf_{\rho_n}, d\tilde {\mj}_n \rangle d\omega  + o(1).
       \eqrs
     
Notice that $\int_{A_{\rho_n}\cap B_{r_n}} -\langle d\mf_{\rho_n}, d\tilde {\mj}_n \rangle d\omega \rightarrow 0$ as $n\rightarrow \infty$ with $r_n\rightarrow 0$, so 
      \bqrs
          \int_{A_{\rho_n}\backslash B_{r_n}}|d(\mh_{\rho_n}-\tilde{\mh}_n)|^2d\omega 
           & = &\int_{A_{\rho_n}}
                   -\langle d\mf_{\rho_n}, d\tilde {\mj}_n \rangle d\omega  + o(1)\\
           & \le & g_1(x^1_n,x^2_n,\rho_n)\|\xi^1_n\|_{\frac{1}{2},2;0} + o(1). 
     \eqrs
Therefore $\lim_{n\rightarrow\infty}\int_{A_{\rho_n}\backslash B_{r_n}}|d(\mh_{\rho_n}-\tilde{\mh}_n)|^2d\omega=0$ and $x^1_n \rightarrow x^1$ strongly in $\hs\cap C^0(\partial B,\mathbb R^k)$. Moreover, by Lemma \ref{tangent} the $N$-harmonic map $\mfr(x^1,a^2)$ is well defined.\medskip 
   
I-c) We shall investigate the behaviour of Jacobi fields.\medskip

For large $n\ge n_0$, $\widetilde{\exp}_{x^1}\eta^1_n=x^1_n$ for some $\eta^1_n\in \mathcal T_{x^1}$, with $\|d\widetilde{\exp}_{x^1,\eta^1_n}\phi^1\|<l_1, \|\phi^1\| < l_1$. Since  $x^1_n \rightarrow x^1$ in $\hs\cap C^0(\partial B,\mathbb R^k)$, $d\widetilde{\exp}_{x^1,\eta^1_n}\rightarrow Id$ in $\hs\cap C^0$. Defining $(v^{\alpha}_n\frac{\partial}{\partial y^{\alpha}}\circ \mf_{\rho_n}) :=\mj_{\mf_{\rho_n}}(d\widetilde{\exp}_{x^1,\eta^1_n}\phi^1,0)$ we have  
           \[\int_{A_{\rho_n}}
               h_{\alpha\beta}\circ \mf_{\rho_n} v_{n,i}^{\alpha}v_{n,i}^{\beta}d\omega 
               \le C \en \text{independent of}\en n\ge n_0.\]
From the Courant-Lebesgue Lemma and $ v^{\alpha}_n|_{\partial_{B_{\rho_n}}}\equiv 0$, 
          \[\int_{\partial(B_{\tilde{r_n}}\cap A_{\rho_n})} 
              h_{\alpha\beta}\circ\mf_{\rho_n}
              \partial_{\theta}v_{n}^{\alpha}\partial_{\theta}v_{n}^{\beta}
              d\theta \le \frac{C}{|\ln \delta|}\en\text{and} 
             \en \|(v^{\alpha}_n)\|_{C^0(B_{\tilde{r_n}(\theta_0)}\cap A_{\rho_n})} 
              \le \frac{C}{|\ln \delta|} \]
for some $\tilde{r_n}\in (\sqrt{\delta}, \sqrt{\sqrt{\delta}})$. Hence, from Lemma \ref{jacobifield}, $E(\mj_{\mf_{\rho_n}}(d\widetilde{\exp}_{x^1,\eta^1_n}\phi^1,0)|_{B_{\tilde{r_n}}})$ is less than $\frac{C}{|\ln \delta|}$. The same holds for $E(\mj_{\mf_{\rho_n}}(d\widetilde{\exp}_{x^1,\eta^1_n}\phi^1,0)|_{B_{r_n}})$, since $r_n\le \tilde{r_n}$. Now choose $\delta$ so that $ \frac{C}{|\ln \delta|} \le \varepsilon:=\frac{1}{n}$.\medskip
 
I-d) Let $\mf_{\rho_n}:=\mf_{\rho_n}(x^1_n,x^2_n)$. The H\"older inequality gives 
        \bqrs
          0 & = & \lim_{n\rightarrow \infty}g^1(x^1_n,x^2_n,\rho_n)\\
           & \ge & \lim_{n\rightarrow \infty}
                    \Big( -\int_{A_{\rho_n}
                            \backslash B_{r_n}}\langle d\mf_{\rho_n}, d\mj_{\mf_{\rho_n}}
                                 (d\widetilde{\exp}_{x^1,\eta^1_n}\phi^1,0)\rangle d\omega 
                          -\int_{B_{r_n}}
                             \langle d\mf_{\rho_n}, d\mj_{\mf_{\rho_n}}
                                   (d\widetilde{\exp}_{x^1,\eta^1_n}\phi^1,0)\rangle d\omega \Big)\\
           & = & \lim_{n\rightarrow \infty}
                   \Big(-\int_{A_{\rho_n}
                           \backslash B_{r_n}}\langle d\mf_{\rho_n}, d\mj_{\mf_{\rho_n}}
                                 (d\widetilde{\exp}_{x^1,\eta^1_n}\phi^1,0)\rangle d\omega - o(1)\Big)\\
          & = & - \int_{A_{\rho}}\langle d\mf, d\mj_{\mf}(\phi^1,0)\rangle d\omega. 
        \eqrs 

The computation in Theorem \ref{appendixthm} yields $ \mf := \mfr(x^1,a_2)\in H^{2,2}(\ar,N)$. Just as in Proposition \ref{h22} we have $\langle\frac{d\mf}{d\vec n},\,\frac{\partial \mf}{\partial\theta}\rangle_h|_{\partial B} \equiv 0$, and clearly $\langle\frac{d\mf}{d\vec n},\,\frac{\partial \mf}{\partial\theta}\rangle_h|_{\partial B_{\rho}}\equiv 0$. \medskip

As a consequence 
      \[\Phi_{\mf}(re^{i\theta}) 
          =r^2\big|\frac{\partial}{\partial r}\mf\big|_h^2
                -\big|\frac{\partial}{\partial \theta} \mf\big|_h^2 
           -2ir\big\langle \frac{\partial}{\partial r}\mf, 
                  \frac{\partial}{\partial \theta}\mf\Big\rangle_h\] 
is real constant.\medskip

Going back to the expression for $\frac{\partial}{\partial \rho} \me$ in Lemma \ref{differ}, the holomorphic function $\Phi_{\mf}$ must be $0$, provided we show $\frac{\partial}{\partial \rho} \me(x^1,a_2,\rho)=0 $. For this one can adapt the argument of \cite{s4}. Thus $\mf := \mfr(x^1,a_2)$ is conformal.\medskip

II) A harmonic, conformal map $\mf := \mfr(x^1,a_2)\in H^{1,2}\cap C^0(\overline{\ar},N)$ must be constant. To prove this we reproduce Theorem 8.2.3 of \cite{j1}.\medskip
  
Consider the complex upper half-plane $ \mathbb C^+ = \{ \theta + ir | r > 0\} $ and let 
             \[\mf((r+\rho)e^{i\theta})=: \widetilde{X}(\theta, r), 
                  \quad \text{well defined on} \en \mathbb R\times [0,1-\rho] \] 
with $\widetilde{X}(\theta,0) = \mf(\rho e^{i\theta}) \equiv a_2$ and $\frac{\partial^m\widetilde{X}}{\partial \theta^{m}} |_{\{r=0\}} \equiv 0$ for each $m$. Choosing an appropriate local coordinate chart in a neighbourhood of $a_2$,  we may assume that $\widetilde{X}(\theta,0) = 0$. Since $\mf$ is conformal and harmonic, $\mf|_{\ar\cup \partial B_{\rho}}\in C^{\infty}$(\cite{hkw}), and by simple computation, $\frac{\partial^m}{\partial \theta^m} \widetilde{X} \equiv \frac{\partial^m}{\partial r^m} \widetilde{X} \equiv 0 \en \text{on} \en{\{r=0\}}, \en m\in \mathbb N$.\medskip

For given $\rho_0\in(0,1)$, let $\Omega := \{ \theta + ir | \theta\in \mathbb R,\en r\in [0, 1-\rho_0)\}$ and $\Omega^- := \{ \theta + ir | \theta\in \mathbb R,\en -r\in [0, 1-\rho_0)\}$. Extending $\widetilde{X}$ to $\Omega\cup\Omega^-=:\widetilde{\Omega}$ by reflection, we see that $\widetilde{X}\in C^{\infty}(\widetilde{\Omega},N)$. From the harmonicity of $\mf$, $|\widetilde{X}_{z\bar z}|\le C|\widetilde{X}_{z}|$ holds. Furthermore $\frac{\partial^m}{\partial \theta^m}\widetilde{X}(0)=\frac{\partial^m}{\partial r^m}\widetilde{X}(0)=0$ and $\lim_{z=(\theta,r)\rightarrow 0} \widetilde{X}(z)|z|^{-m} = 0$ for all $m\in \mathbb N$. Hence $\widetilde{X}$ is constant in $\widetilde{\Omega}$ by the Hartman-Wintmer Lemma (see \cite{j1}). This holds for each $\rho_0\in(0,1)$, so we get $\mf\equiv a_2$ on $\overline{\ar}$. But this contradicts the assumption $\text{dist}(\Gamma_1,\Gamma_2) > 0$. Therefore case 2 cannot really occur.\medskip

{\bf (case3)} Suppose that  $x^i_n=\gamma_i\circ w^i_n \rightarrow$ const.$=:a_i\in \Gamma_i$ in $L^1(\partial B,\mathbb R^k), i=1,2$. Similarly to case 2, this will lead to a contradiction.\medskip

First of all $\Phi_{\mf}$ is a real constant for $\mf:=\mf(a^1,a^2)$. Supposing that $\frac{d}{d\rho}E(\mf)\not= 0$, we have $\big |\int^{2\pi}_{0}\int^{1-t}_{\rho+\delta}\big[\big|\frac{\partial}{\partial r}\mf_{\rho_n}\big|_{h}^2 -\frac{1}{r^2}\big|\frac{\partial}{\partial r}\mf_{\rho_n}\big|_{h}^2\big] \frac{1}{1-\rho-\delta}drd\theta\big |=C>0$ for some fixed $t,\delta >0$ and large $n \ge n_0$. Let 
       \[\widetilde{\mf}_{n}^{\sigma} := \left\{\begin{array}{r@{\quad\quad}l}
                   \mf_{\rho_n} & \text{on}\en A _{1-t}, \\
                   \mf_{\rho_n}\circ \tau^{\rho +\delta}_{\sigma;1-t} & \text{on}\en A_{\sigma}\backslash A _{1-t},\\
                   \mf_{\rho_n}(\frac{\rho+\delta}{\sigma}r,\theta) 
                             & \text{on}\en A_{\frac{\sigma\rho_n}{\rho+\sigma}} \backslash A_{\sigma},
                        \end{array}\right. \]
where $\tau^{\rho +\delta}_{\sigma;1-t}$ is a diffeomorphism from $[\sigma, 1-t]$ to $[\rho + \delta, 1-t]$. Then  
       \[2\frac{d}{d\sigma}E(\widetilde{\mf}_{n}^{\sigma})|_{\sigma = \rho+\delta}
                      =  \int^{2\pi}_{0}\int^{1-t}_{\rho+\delta}
                    \left[|\partial_r \mf_{\rho_n}|^2 - \frac{1}{r^2}|\partial_{\theta}\mf_{\rho_n}|^2\right]
                        \frac{1-t}{1-t-\rho-\delta}drd\theta.\] \\
     Since $\widetilde{\mf}_{n}^{\rho+\delta} = \mf_{\rho_n}$ it follows that         
    \[\rho_n |g_3(x_n) |  =  |\rho_n\frac{d}{d\sigma}E(\mf_{\rho_n})|_{\sigma = \rho_n}|
              =|(\rho+\delta)\frac{d}{d\sigma}E(\widetilde{\mf}_{n}^{\rho+\delta})| \ge C  >0,\]
contradicting the assumption  $g_3(x_n)\rightarrow 0$. Thus, $\mfr(a_1,a_2)$ is conformal, and we can use the argument of {\bf (case2)}-II). \medskip      

{\bf (case4)}: Suppose that  $\rho=0$.\medskip

For conformal diffeomorphisms $\tau^i_n$ of $B$, $\mf^i(x^i_n)\circ\tau^i_n = \mf^i(\widetilde{x^i_n})$ holds with $\widetilde{x^i_n}\in M^{i\ast}$, $i=1,2$. Furthermore $\widetilde{x^i_n}$ has a subquence converging to $x^i \in M^{i\ast}$ uniformly.\medskip

For given $\varepsilon>0$, there exist $\delta>0$ and $n_0$ such that for $n\ge n_0$, $(x^1_n,x^2_n,\rho_n)\in \mathcal N_{\delta}(\widetilde{x^1_n},\widetilde{x^2_n},0)$ and $|g(x^1_n,x^2_n,\rho_n)-g(\widetilde{x^1_n},\widetilde{x^2_n},0)| <\varepsilon$. Thus, from the topology of $\overline{\mm^{\ast}}$ we have then $g(\widetilde{x^1_n},\widetilde{x^2_n},0)\rightarrow 0$ when $n\rightarrow \infty$. \medskip

Once again as in (case 1), some subsequence of $\widetilde{x^i_n}$ strongly converges to $x^i\in M^{i\ast}$ with  $g(x^1,x^2,0)=0$. This finishes the proof of the PS condition.                           \hfill $\Box$


\subsection{Unstable minimal surfaces of annulus type}
This section contains three Lemmata, adapted from \cite{s4} to our purposes, in preparation to the main theorems. Before that, we need some explanation for $\overline{\mm^{\ast}}$ (see section \ref{Palais}).\medskip

I) The boundary $\partial \mm^{\ast}$:
\vspace*{-0.1cm} 
    \begin{enumerate}
    \renewcommand{\labelenumi}{(\roman{enumi})}
      \item
      For an element $x^i\in M^i$, $(x^i)^{-1}(P^i_k)$ is a closed interval on the unit circle, since $x^i$ is weakly monotone. Let $Q^i_k$ be the first endpoint of $(x^i)^{-1}(P^i_k)$ relative to the positive orientation of the circle for each $i=1,2, \en k=1,2,3$. Taking the conformal linear fractional transformation $T_{x^i}$ of the unit disc which maps $(\cos \frac{2\pi (k-1)}{3}, \sin \frac{2\pi (k-1)}{3})$ to $Q^i_k$ and the unit circle onto itself, we have $x^i\circ T_{x^i}\in M^{i\ast}$. Moreover $T_{x^i\circ T_{x^i}}=Id$, since $T_{x^i}$ is one-to-one.\smallskip

      For $x^i, y^i \in M^i$, we write $x^i\sim y^i$ if $x^i\circ T_{x^i} = y^i\circ T_{y^i}$, clearly an equivalence relation.  Now we can quotient $M^i$ in such a way that each class possesses exactly one element from 
      $x^i\in M^{i\ast}$, denoted by $[x^i]\in M^{i\ast}$, with $\|[x^i]\| = \|x^i\|$.

       \item
       For $x^i\in M^{i\ast}$ and $\xi^i\in\mathcal T_{x^i}, \|\xi^i\| < l_i$ we may calculate $[\widetilde{\exp}]_{[x^i]}\xi^i := [\widetilde{\exp}_{x^i}\xi^i] = [\tilde{x}^i]\in M^{i\ast}$, where  $\widetilde{\exp}_{x^i}\xi^i\circ T_{\widetilde{\exp}_{x^i}\xi^i} = \tilde{x}^i\in M^{i\ast}$. We will denote this correspondence simply by $\widetilde{\exp}_{x^i}\xi^i = \tilde{x}^i$, which is clearly smooth, since $T_{x^i}$ varies smoothly with $x^i \in M^i$ (cf. above).\smallskip

       Now, for $[x]=([x^1],[x^2],0) \in \partial \mm^{\ast}$ with $x^i\in M^{i\ast}$, we define $g([x]):= g(x)$, where $ x=(x^1, x^2, 0)$. Recall that the Dirichlet integral is invariant under conformal mappings, so for $\xi^1\in \mathcal T_{x^1}$
          \bqrs
          \lefteqn{ \me([\widetilde{\exp}]_{[x^1]}t\xi^1, [x^2],0)}\hspace*{1.0cm}\\ 
             \hspace*{0.5cm} &=& E(\mf^1([\widetilde{\exp}]_{[x^1]}t\xi^1)) + E(\mf^2([x^2]))
                                            = E(\mf^1([\widetilde{\exp}_{x^1}t\xi^1])) + E(\mf^2(x^2))\\
                           &=& E(\mf^1(\tilde{x}^1_t))+ E(\mf^2(x^2)) 
                                           = E(\mf^1(\widetilde{\exp}_{x^1}t\xi^1)) + E(\mf^2(x^2))\\
                      &=& \me(\widetilde{\exp}_{x^1}t\xi^1, x^2,0),
          \eqrs

 
        where $\widetilde{\exp}_{x^1}t\xi^1\circ T_{\widetilde{\exp}_{x^1}t\xi^1} = \tilde{x}^1_t\in M^{1\ast}$. The same holds for $\me(x^1, \widetilde{\exp}_{x^2}t\xi^2,0), \en \xi^2\in \mathcal T_{x^2}$. Therefore, $g([x])$ is well defined. 
    \end{enumerate}

II) The interior $\mm^{\ast}$:
    \begin{enumerate}
    \renewcommand{\labelenumi}{(\roman{enumi})} 
        \item
           For $x=(x^1, x^2, \rho)\in \mm$ let $Q^1_1$ be the first endpoint of $(x^1)^{-1}(P^1_1)$ relative to the positive orientation on the circle, and $r_{x^1}$ the positive rotation of $A_{\rho}$ mapping the point $(1,0)$ to $Q^1_1$ of the unit circle. Then $x\circ r_{x^1}:=(x^1\circ r_{x^1}, x^2\circ r_{x^1}, \rho)\in \mm ^{\ast}$ for each $\rho\in (0,1)$ and $r_{x^1\circ r_{x^1}}=Id$.\smallskip

           Since $(x^1, x^2, \rho)$ and $(x^1\circ r_{x^1}, x^2\circ r_{x^1}, \rho)$ can be naturally identified, it makes sence to define an equivalence relation $x\sim y$ if $x\circ r_{x^1}= y\circ r_{y^1}$ holds, for $x, y\in \mm$. In each equivalence class there is exactly one element from $\mm^{\ast}$. 
       \item
          For $[x] = ([x^1, x^2],\rho)$ with $(x^1, x^2, \rho)\in \mm^{\ast}$ and $\xi^1\in \mathcal T_{x^1}$, we compute $[\widetilde{\exp}]_{[x]}(\xi^1,0,0) =([\widetilde{\exp}_{x^1}\xi^1,x^2],\rho)=([\widetilde{\exp}_{x^1}\xi^1\circ r,x^2\circ r],\rho)$, with $r:=r_{\widetilde{\exp}_{x^1}\xi^1}$. Denoting this correspondence simply by $\widetilde{\exp}_{x}(\xi^1,0,0)=(\widetilde{\exp}_{x^1}\xi^1\circ r,x^2\circ r,\rho)\in \mm^{\ast}$, $\widetilde{\exp}$ is clearly smooth.\smallskip
 
          Let $g([x]):=g(x)$ for $x\in \mm^{\ast}$ and observe that for $\xi^1\in \mathcal T_{x^1}$ 
           \bqrs
             \me([\widetilde{\exp}]_{[x]}(t\xi^1,0,0))  &=&\me([\widetilde{\exp}_{x^1}t\xi^1, x^2], \rho)
                          = E(\mfr([\widetilde{\exp}_{x^1}t\xi^1,x^2] ))\\
                      &=& E(\mfr( \widetilde{\exp}_{x^1}t\xi^1\circ  r_t, x^2 \circ  r_t)
                          = E((\mfr(\widetilde{\exp}_{x^1}t\xi^1, x^2))\circ  r_t)\\ 
                      &=& E(\mfr(\widetilde{\exp}_{x^1}t\xi^1,x^2 )) = \me(\widetilde{\exp}_{x^1}t\xi^1, x^2,\rho),  
          \eqrs
        where $r_t:=r_{\widetilde{\exp}_{x^1}t\xi^1}$ is such that $(\widetilde{\exp}_{x^1}t\xi^1\circ  r_t, x^2 \circ  r_t, \rho)\in \mm^{\ast}$. Moreover, for $\xi^2\in \mathcal T_{x^2}$, $[\widetilde{\exp}]_{[x]}(0,\xi^2,0) = ([x^1, \widetilde{\exp}_{x^2}\xi^2],\rho)$ with $(x^1, \widetilde{\exp}_{x^2}\xi^2,\rho)\in\mm^{\ast}$, so we can compute as usual, and $g([x])$ is well defined.
      \end{enumerate}

\begin{remark}\label{hihi}
Consider $x\in\partial\mm ^{\ast}$ and $y\in\mm ^{\ast}$ as equivalence classes. Then for $\xi = (\xi^1, \xi^2, 0)\in\mathcal T_{x}\partial\mm$ and $\xi_{\rho} = (\xi^1, \xi^2, \rho)\in\mathcal T_{x}\mm$ with $\|\xi^i\|_{\hfn} \le l_i$, we have $\widetilde{\exp}_x \xi \in\partial\mm ^{\ast}$ and $\widetilde{\exp}_y \xi_{\rho} \in\mm ^{\ast}$. Moreover, $\widetilde{\exp}_x $ and $d\widetilde{\exp}_x$ are continuous.
 \end{remark}

And now we proceed with the results.

\begin{lemma}\label{vector}
For any $\delta > 0$, there exists a uniformly bounded, continuous vector field $e_{\delta}:M^1\times  M^2\times [0,1)\rightarrow \mathcal T_{M^1}\times \mathcal T_{M^2}\times \mathbb R$, satisfying local Lipschitz continuity on $\mm$ and $\partial \mm$(separably) with the following properties
\begin{enumerate}
  \renewcommand{\labelenumi}{(\roman{enumi})}
   \item
     for $\beta\in \mathbb R$ there exists $\varepsilon > 0$ such that $y_{\delta}(x) = \left(\widetilde{\exp}_{x^1}e_{\delta}^1(x^1),\widetilde{\exp}_{x^2}e_{\delta}^2(x^2), \rho +e_{\delta}^3(\rho)\right) \in \mm(\rho):= \{ x\in \mm | x=(x^1,x^2,\rho) \}$ for any $x\in \mm(\rho)$ with $\me(x) \le \beta$ and $0<\rho<\varepsilon$ (that is, $e_{\delta}$ is parallel to $\partial\mm$ near $\partial \mm$),
   \item
     for any such $\beta, \me, x$ and any pair $T=(\tau^1,\tau^2)$ of conformal transformations of $B$,  $ y_{\delta}(x\circ T) = y_{\delta}(x)\circ T$, where $ \mf^i\left((x\circ T)^i\right ) = \mf^i(x^i)\circ T, \quad i=1,2$,
   \item
     for any $x\in \overline {\mm}$, $ \langle d\me(x),e_{\delta}(x)\rangle_{\mathcal T_{x^1}\times \mathcal T_{x^2} \times \mathbb R} \le \delta-g(x)$,
   \item
     for $x\in \mm^{\ast}$ and $y\in\partial\mm^{\ast}$, we have $y_{\delta}(x)\in \mm^{\ast}$ and $y_{\delta}(y)\in\partial\mm^{\ast}$. 
  \end{enumerate}       
\end{lemma}
{\bf Proof.} The proof of the analogous result in \cite{s4} can easily be adapted to our setting, because Remark \ref{hihi} holds.

\begin{lemma}\label{flow}
For a given locally Lipschitz continuous vector field $f:\overline{\mm} \rightarrow \mathcal T_{M^1}\times \mathcal T_{M^2}\times \mathbb R$ satisfying Lemma \ref{vector}, there exists a unique flow $\Phi : [0,\infty)\times \overline {\mm^{\ast}} \rightarrow \overline{\mm^{\ast}}$ with 
  \[\Phi(0,x) = x,\en \frac{\partial}{\partial t}\Phi(t,x)=f\left(\Phi(t,x)\right),\quad x\in 
         \overline {\mm^{\ast}}\,.\] 
\end{lemma}
{\bf Proof.} We use Euler's method. Define $\Phi^{(m)}:[0,\infty)\times \overline {\mm^{\ast}} \rightarrow \overline{\mm^{\ast}},\en m\ge m_0$, by  
      \bqr
          \Phi^{(m)}(0,x) &:= & x\nonumber \\
          \Phi^{(m)}(t,x) & :=&\widetilde{\exp}_{\Phi^{(m)}(\frac{[mt]}{m},x)}
             \Big(\frac{mt-[mt]}{m} f\big(\Phi^{(m)}(\frac{[mt]}{m},x)\big)\Big), \en t>0 \label{tra}
      \eqr   
where $[\tau]$ denotes the largest integer which is smaller than $\tau\in \mathbb R$. This is well defined due to the convexity of $\mathcal T_{x^i},\, x^i\in M^i, i=1,2$ and by Lemma \ref{vector} {\it(iv)}.\medskip
  
Recalling the map $w^i\in C^0(\mathbb R, \mathbb R)$ with $x^i = \gamma^i\circ w^i ,\en x^i\in M^i$ (section \ref{construction} {\bf III)}), consider 
      \[ W^i:=\{w^i\in C^0(\mathbb R,\mathbb R):\gamma^i\circ w^i=x^i\en 
         \text{for some}\en x^i\in M^i\}, \quad W:= W^1\times W^2\times [0,\infty). \]
Let $ \gamma(w):= (\gamma^1\circ w^1,\gamma^2\circ w^2,\rho) \en\text{for}\en (w^1,w^2, \rho)=: w \in W, \en \gamma :=(\gamma^1, \gamma^2, Id)$ and $\widetilde{f} := (\widetilde{f}^1,\widetilde{f}^2, f^3) \en\text{with}\en \widetilde{f^i}(w^i):=(d\gamma^i)^{-1}(f^i(x^i))\in C^0(\mathbb R / 2\pi,\mathbb R)$. Then there exists $\widetilde{\Phi}^{(m)}(t,w)\in W$ with $\Phi^{(m)}(t,x)=\gamma(\widetilde{\Phi}^{(m)}(t,w))$, so we can write (\ref{tra}) as follows: 
     \bqrs
        \widetilde{\Phi}^{(m)}(t,w) &=& \widetilde{\Phi}^{(m)}(\frac{[mt]}{m},w) + 
             \frac{mt-[mt]}{m} \widetilde{f}\big(\widetilde{\Phi}^{(m)}(\frac{[mt]}{m},w)\big) + 2\pi l, \en l\in Z.
    \eqrs

When $t\in(\frac{k}{m}, \frac{k+1}{m}], k\in \mathbb Z$, $\widetilde{\Phi}^{(m)}(t,w) = \widetilde{\Phi}^{(m)}(0,w) + 
 \int^{t}_{0}\widetilde f(\widetilde{\Phi}^{(m)}(\frac{[ms]}{m},w))ds$.\medskip

The compution now proceeds as in the Euclidean case. For any $T>0$, $G>0$, there exists $C(T,G)$ with $\|\Phi^{(m)}(\cdot,w)\|_{L^{\infty}([0,T]\times W_K,\overline{\mm})}  \le C(T,G), \en w\in W \en \text{with}\en \|w\|_{W}\le G$.\medskip

Let $L_1$ resp. $L_2$ be the Lipschitz constants of $f$ in $\{x\in \mm \,|\,\|x\|\le C(T,G)\}$ and $\{ x\in \partial\mm \,|\,\|x\|\le C(T,G)\}$, and call $L:= \max\{C(\gamma^i)L_1, C(\gamma^i)L_2\}$.\medskip 
 
For $\frac{m}{n} < 1$, $\en \|\widetilde{\Phi}^{(m)}(t,w)-\widetilde{\Phi}^{(n)}(t,w)\| \le tL\frac{2}{m}C(f)+ tL\|\widetilde{\Phi} ^{(m)}(\cdot,w)- \widetilde{\Phi} ^{(n)}(\cdot,w)\|_{L^{\infty}([0,t],W)}.$ Hence, for $m,n\ge m_0$, 
we have \medskip

$\|\widetilde{\Phi}^{(m)}(\cdot,w)-\widetilde{\Phi}^{(n)}(\cdot,w)\|_{L^{\infty}([0,t],W)} \le tL(\frac{2}{m}+\frac{2}{n})C(f) + tL\| \widetilde{\Phi}^{(m)}(\cdot,w) - \widetilde{\Phi}^{(n)}(\cdot,w)\|_{L^{\infty}([0,t],W)}$.\medskip 

By choosing $t\le \min\{T,\frac{1}{2L}\}$, $\{\widetilde{\Phi}^{(m)}\}$ converges uniformly to some function $\widetilde{\Phi}$ on $[0,t]\times \{w\in W: \|w\| \le G\}$ as $m\rightarrow \infty$. Then $\frac{\partial}{\partial t}\widetilde{\Phi}(t,w) = \widetilde f\big(\widetilde{\Phi}(t,w)\big)$. \medskip

For $\Phi(t,w):= \gamma\circ\widetilde{\Phi}(t,w)\in \overline{\mm}^{\ast}$, the uniform boundedness of $f$ yields a flow $\Phi$ such that $\frac{\partial}{\partial t}{\Phi}(t,w) = d\gamma\Big(\widetilde f\big(\widetilde{\Phi}^{(m)}(t,w)\big)\Big)= f\big({\Phi}(t,w)\big)$ for each $x\in \overline{\mm}$. $\Phi(t,w)$ depends continuously on the initial data, and it can be prolonged for $t>0$.                                          \hfill $\Box$ \medskip

The next result is slightly weaker than the corresponding Lemma 4.15 in \cite{s4}, but will nevertheless suffice for our aim. 

 \begin{lemma}\label{boundary}
   Let $\mf^i(x^i_0)$ be a solution of $\mathcal P(\Gamma_i)$  for some $x_0^i\in M^i$, $i=1,2$, and suppose that $d:=\text{dist}(\mf^1(x^1_0),\mf^2(x^2_0)) > 0$. Then there exist 
   $\varepsilon >0$, $\rho_0\in(0,1)$ and $C>0$, dependent on $\me(x^1_0,x^2_0,0)$ 
   such that for $x^i\in M^i$ with $\|x^i - x^i_0\|_{\frac{1}{2},2;0} =:s(x^i) < \varepsilon$,
   \[\me(x^1,x^2,\rho) \ge \me(x^1,x^2,0)+\frac{Cd^2}{|\ln \rho|}, \quad \text{for all}\en \rho\in (0,\rho_0).\]
 \end{lemma}
{\bf Proof.} Let $\mfr:=\mfr(x^1,x^2)$ $\mf^i:=\mf^i(x^i), i=1,2$.  Choose $\sigma_1$ and $\delta$ such that $\sqrt{\rho}<\delta<\sigma_1<\sqrt{\sqrt{\rho}}$. For $T(re^{i\theta}):=\rho\frac{1}{re^{i\theta}}$ and $\sigma_2 := \frac{\rho}{\delta}$, take $f_{\sigma_1} := \mfr|_{A_{\sigma_1}}$ and $g_{\sigma_2} := \mfr|_{B_{\delta}\backslash B_{\rho}}(T^{-1})$. Then 
    \begin{equation}\label{hoohoo}
           E(\mfr)  =  E(f_{\sigma_1})+E(\mfr|_{B_{\sigma_1}\backslash B_{\delta}})+E(g_{\sigma_2}).
    \end{equation}
     
We will estimate $E(\mfr)$ in (I) {\bf ---} (III).\medskip

(I) Estimate of $E(f_{\sigma_1})$ and $E(g_{\sigma_2})$.\medskip

In order to control $E(f_{\sigma_1})$ we take $a_1\in N$ with $\min_{a\in N} E(\mf_{\sigma_1}(x^1,a)) = E(\mf_{\sigma_1}(x^1,a_1))$ and let $\mf^1_{\sigma_1}:=\mf_{\sigma_1}(x^1,a_1)$.\medskip

Next, define $\widetilde{\mf^1_{\sigma_1}}:B\rightarrow N$ as follows: 
Let $\widetilde{\mf^1_{\sigma_1}}|_{B\backslash B_{\frac{1}{2}}}$ be $\mf^1_{B\backslash B_{\frac{1}{2}}}$, $\widetilde{\mf^1_{\sigma_1}}|_{B_{\frac{1}{2}}\backslash B_{\sigma_1}}$ be harmonic on $N$ with $\mf^1|_{\partial B_{\frac{1}{2}}}$ on $\partial B_{\frac{1}{2}}$ and $\mf^1(0)$ on $\partial B_{\sigma_1}$, and set  $\widetilde{\mf^1_{\sigma_1}}|_{B_{\sigma_1}} \equiv \mf^1(0)$. Thus 
       \[ 2E(\widetilde{\mf^1_{\sigma_1}}-\mf^1)
         =\underbrace{\int_{B_{\frac{1}{2}}\backslash B_{\sigma_1}}|\nabla
          (\widetilde{\mf^1_{\sigma_1}}-\mf^1)|^2d\omega}_{=:a} 
         + \underbrace{\int_{B_{\sigma_1}}|\nabla (\widetilde{\mf^1_{\sigma_1}}-\mf^1)|^2d\omega}_{=:b}.\]
It is easy to see that $b\le C|\sigma_1|^2$, since $\mf^1$ is regular on $B_{\frac{1}{2}}$. Notice $\widetilde{\mf^1_{\sigma_1}}|_{B_{\frac{1}{2}}\backslash B_{\sigma_1}}\in H^{2,2}$, since $\widetilde{\mf^1_{\sigma_1}}|_{\partial B_{\frac{1}{2}}}$ is regular and constant on $\partial B_{\sigma_1}$. Thus, 
    \bqrs 
     a & = & \int_{\partial B_{\frac{1}{2}}}\langle \nabla (\widetilde{\mf^1_{\sigma_1}}-\mf^1)\vec n,
               \widetilde{\mf^1_{\sigma_1}}-\mf^1\rangle d_0 +\int_{\partial B_{\sigma_1}}
          \langle \nabla (\widetilde{\mf^1_{\sigma_1}}-\mf^1)\vec n,\widetilde{\mf^1_{\sigma_1}}-\mf^1\rangle d_0\\
      &\le& C\|\mf^1(0)-\mf^1|_{\partial B_{\sigma_1}}\|_{C^0}\sigma_1 \le C|\sigma_1|^2,
           \quad \text{with}\en   C=C(E(\mf^1(x^1))).
    \eqrs 
Let $\mf^1_{\sigma_1}|_{B_{\sigma_1}}\equiv a_1$,  so that $E(\mf^1)\le E(\mf^1_{\sigma_1}) \le E(\widetilde{\mf^1_{\sigma_1}})$. From Lemma \ref{huhu!}, 
        \bqr
           \lefteqn{E(\mf^1_{\sigma_1}- \mf^1)\le E(\mf^1_{\sigma_1}) - E(\mf^1) + o_s(1) }
                  \hspace{0.3cm}\label{vorsicht}\\
                 & &  \le E(\widetilde{\mf^1_{\sigma_1}}) - E(\mf^1)+o_s(1) 
                     \le E(\widetilde{\mf^1_{\sigma_1}}-\mf^1) + o_s(1)\le C|\sigma_1|^2 + o_s(1)\nonumber,
        \eqr
where $o_s(1)\rightarrow 0$ as $\|x^1 - x^1_0\|_{\frac{1}{2},2;0} =:s(x^1)\rightarrow 0$.\medskip

Since $E(\mf^1_{\sigma_1}- \mf^1)|_{B_{\sigma_1}} \le C|\sigma_1|^2+ o_s(1)$, we have $E(\mf^1_{\sigma_1}- \mf^1)|_{A_{\sigma_1}} \le C|\sigma_1|^2+ o_s(1)$.\medskip
        
For $X^1 := f_{\sigma_1} -  \mf^1_{\sigma_1}$,  
    \[\left|\int_{A_{\sigma_1}}\nabla (\mf^1_{\sigma_1}- \mf^1)\nabla X^1 d\omega \right| 
      \le  C\sigma_1\left(\int_{A_{\sigma_1}}\left|\nabla(f_{\sigma_1} - \mf^1_{\sigma_1})\right|^2 d\omega \right)^{\frac{1}{2}} \le  C\sigma_1.\]\medskip
         
On the other hand, 
       \bqr
         \lefteqn{|a_1-\mf^1(0)|^2 = 
          \left|\int^1_{\sigma_1}\partial_r(\widetilde{\mf^1_{\sigma_1}}-\mf^1_{\sigma_1})dr\right|^2
          \le(1-\sigma_1)\int^1_{\sigma_1}\left|\nabla(\widetilde{\mf^1_{\sigma_1}}-\mf^1_{\sigma_1})\right|^2 dr}
          \nonumber \hspace{3.5cm}\\
            &\le&\frac{1-\sigma_1}{\sigma_1}(E(\widetilde{\mf^1_{\sigma_1}}- \mf^1) +E(\mf^1_{\sigma_1}- \mf^1)) 
                         \le C \sigma_1\label{hoo!} + o_{s}(1).
       \eqr
          
From the above consideratoins  
         \bqrs
          \Big|\int_{A_{\sigma_1}}\langle\nabla \mf^1_{\sigma_1},\nabla X^1 \rangle d\omega \Big| 
             & \le &  \Big|\int_{A_{\sigma_1}}\langle\nabla \mf^1,\nabla X^1 \rangle d\omega \Big| + C\sigma_1 \\
             & \le &  \|\nabla \mf^1|_{\partial B_{\sigma_1}}\|\|(-a_1+\mfr|_{\partial B_{\sigma_1}})\|
                                         \sigma_1+C\sigma_1  \le  C\sigma_1.
         \eqrs
        
With $C\in \mathbb R$ depending on $E(\mf^1)$ 
           \bqr
              E(f_{\sigma_1})  
                =E(\mf^1_{\sigma_1}) +  \int_{A_{\sigma_1}}\langle\nabla \mf^1_{\sigma_1},
                     \nabla X^1 \rangle d\omega + E(X^1)
               \ge E(\mf^1)  - C\sigma_1.\label{estimate1} 
           \eqr  
Similarly $E(g_{\sigma_2}) \ge  E(\mf^2)  - C\sigma_2$, and $C$ depends on $E(\mf^2)$.\bigskip

(II) Estimate of $E(\mfr|_{B_{\sigma_1}\backslash B_{\delta}})$.\medskip 
                            
From (\ref{hoo!}), $|a_1-a_2| \ge \left||\mf^1(0)-\mf^2(0)|-|a_1-\mf^1(0) + \mf^2(0)-a_2|\right|\ge d - o_{\rho}(1) - o_s(1)$.\medskip
           
Let $\mathcal H^a_b(f,g)$ be the harmonic map on $B_{a}\backslash B_{b}$ in $\mathbb R^k$ with boundary $f$ on $\partial B_{a}$ and $g$ on  $\partial B_{b}$.\medskip

Writing $\sigma_1=:\sigma,\en \frac{\delta}{\sigma^1}=:\tau,\en \mfr|_{\partial B_{\sigma_1}}=:p, \en \mfr|_{\partial B_{\delta}}=:q$, we have 
        \bqrs
           \lefteqn{\Big|\int \langle \nabla \mathcal H^{\sigma}_{\delta}(a_1,a_2), 
                      \nabla \mathcal H^{\sigma}_{\delta}(-a_1+p,-a_2+q)\rangle d\omega \Big|}\\
             &= &  \Big|\int \langle \nabla \mathcal H^{1}_{\tau}(0,-a_1+a_2), 
                     \nabla \mathcal H^{1}_{\tau}(-a_1+p(\cdot \sigma),-a_2+q(\cdot\sigma)\rangle d\omega \Big|\\
             &\le& \frac{2\pi}{|\ln\tau|}|-a_1+a_2|\big(|-a_1+p(\cdot \sigma)|+|-a_2+q(\cdot \sigma)|\big)  
                   \le  C\frac{(o_{\rho}(1)+o_s(1))}{|\ln\rho|}.
        \eqrs
   
Moreover 
        \bqrs
          E( \mathcal H^{\sigma}_{\delta}(a_1,a_2)) &\ge&  
          E(\mathcal H^{1}_{\rho}(0, -a_1+a_2)) = E\big(( -a_1+a_2)\frac{\ln r}{\ln \rho}\big)
                  \ge \frac{\pi d^2}{|\ln \rho|}-C\frac{(o_{\rho}(1)+o_s(1))}{|\ln \rho|}.
        \eqrs     
Thus, 
           \bqr
             E(\mfr|_{B_{\sigma}\backslash B_{\delta}}) 
                 &\ge& E(\mathcal H^{\sigma}_{\delta}(p,q))
                        =  E\big(\mathcal H^{\sigma}_{\delta}(a_1,a_2) + \mathcal 
                         H^{\sigma}_{\delta}(-a_1 + p,-a_2 + q)\big) \nonumber\\
                 & \ge & \frac{\pi d^2}{|\ln \rho|}-C\frac{o_{\rho}(1)+o_s(1)}{|\ln \rho|}\label{estimate2}
            \eqr        
with $C$ depending only on  $E(\mf^i)$, $i=1,2$.\medskip

(III) Estimate $E(\mfr)$.\medskip

From (\ref{hoohoo}), (\ref{estimate1}), (\ref{estimate2}) and the choice made for $\sigma_i$, $i=1,2$,   
           \bqrs 
              \lefteqn{\me(x^1,x^2,\rho) \ge \me(x^1,x^2,0)-C\sigma_i + \frac{\pi d^2}{|\ln \rho|}-
             C\frac{(o_{\rho(1)}+o_{s(1)})}{|\ln \rho|}}\hspace{0.5cm}\\
                                &\ge& \me(x^1,x^2,0)-C(\sqrt{\rho}+\sqrt{\sqrt{\rho}}) + \frac{\pi d^2}{|\ln \rho|}-
                                      C\frac{(o_{\rho}(1)+o_{s(1)})}{|\ln \rho|}
                                 \ge \me(x^1,x^2,0) + C\frac{d^2}{|\ln \rho|},
           \eqrs
for $\rho\le \rho_0$, some small $\rho_0\in(0,1)$ and $s(x^i)$.\hfill $\Box$

\begin{lemma}\label{huhu!}
  With the same notations as in  Lemma \ref{boundary}, 
    \[E(\mf^1_{\sigma_1}- \mf^1)= E(\mf^1_{\sigma_1}) - E(\mf^1) + o_s(1).\] 
\end{lemma}
{\bf Proof.} Let $G^1:=\mf^1(x^1_0)$. Note $\min_{a\in N} E(\mf_{\sigma_1}(x^1_0,a))=E(\mf_{\sigma_1}(x^1_0,a^1))$, and let $G^1_{\sigma_1}:=\mf_{\sigma_1}(x^1_0,a^1)$,  $G^1_{\sigma_1}|_{B_{\sigma_1}}\equiv a^1$.\medskip

Since $G^1\in H^{2,2}$ 
           \[0= \int_B \langle \nabla G^1, \nabla(G^1_{\sigma_1}- G^1)\rangle d\omega =
                        \int_B \langle II\circ G^1(dG^1,dG^1),G^1_{\sigma_1}- G^1\rangle.\]
Note that $\|\mf^1_{\sigma_1} - G^1_{\sigma_1}\|_{C^0} \rightarrow 0$ when $\|x_0^1 - x^1\|_{\frac{1}{2},2;0}=:s(x^1) \rightarrow 0$ just as in Lemma \ref{differ} {\it(B)}. Moreover, $\|G^1 - \mf^1\|_{1,2;0}\rightarrow 0$ as $s(x^1) \rightarrow 0$, so by the H\"older inequality,  
          \[\left|\int_B \langle II\circ\mf^1(d\mf^1,d\mf^1),\mf^1_{\sigma_1}- \mf^1\rangle d\omega 
             - \int_B \langle II\circ G^1(dG^1,dG^1),G^1_{\sigma_1}- G^1\rangle d\omega\right| = o_s(1).\]

In this way 
         \[ \int_B \langle \nabla \mf^1, \nabla(\mf^1_{\sigma_1}- \mf^1)\rangle d\omega =
            \int_B \langle II\circ\mf^1(d\mf^1,d\mf^1),\mf^1_{\sigma_1}- \mf^1\rangle = o_s(1)\]
and 
          \bqrs
            2E(\mf^1_{\sigma_1}- \mf^1)
                 &  = & \int_B \langle \nabla \mf^1_{\sigma_1}, \nabla(\mf^1_{\sigma_1}- \mf^1)\rangle d\omega
                                -\int_B \langle \nabla \mf^1, \nabla(\mf^1_{\sigma_1}- \mf^1)\rangle d\omega \\
                & = & \int_B \langle \nabla \mf^1_{\sigma_1}, \nabla(\mf^1_{\sigma_1}- \mf^1)\rangle d\omega 
                                 + o_s(1)  \\
                   & = & \int_B |\nabla \mf^1_{\sigma_1}|^2d\omega - 
                    \int_B\langle \nabla \mf^1,\nabla \mf^1_{\sigma_1}-\nabla\mf^1\rangle d\omega 
                               - \int_B|\nabla \mf^1|^2d\omega   + o_s(1) \\ 
                  & = &\int_B |\nabla \mf^1_{\sigma_1}|^2d\omega - \int_B|\nabla \mf^1|^2d\omega   + o_s(1).
          \eqrs 
                                                                                     \hfill $\Box$\medskip
        
We eventually arrive at 

\begin{theorem}  
Let $\Gamma_1, \Gamma_2 \subset (N,h)$ satisfy (C1) or (C2) and define  
       \bqrs
               d  &=&  \inf\{ E(X)\, | X\in \mathcal S (\Gamma_1,\Gamma_2)\} \\
               d^{\ast} &=&  \inf \{ E(X^1) + E(X^2)\, | X^i \in  \mathcal S (\Gamma_i), i=1,2 \}.
       \eqrs
If $d<d^{\ast}$, there exists a minimal surface of annulus type bounded by $\Gamma_1$ and $\Gamma_2$.
\end{theorem}
{\bf Proof.}  The PS condition (Proposition \ref{psc}) and Proposition \ref{h22} allow to conclude straight away. For details we refer to \cite{s1}.  \hfill $\Box$

\begin{theorem} 
Let $\mf^1$, resp. $\mf^2$, be an absolute minimizer of $E$ in $\mathcal S(\Gamma_1)$, resp. $\mathcal S(\Gamma_2)$, and suppose that $\text{dist}(\mf^1,\mf^2) > 0$. Assume furthermore there is a strict relative minimizer of $E$ in $\mathcal S(\Gamma_1, \Gamma_2)$. Then there exists either a solution of $\mathcal P(\Gamma_1,\Gamma_2)$ which is not a relative minimizer of $E$ in $\mathcal S(\Gamma_1, \Gamma_2)$, i.e. an unstable annulus-type-minimal surface, or a pair of solutions to $\mathcal P(\Gamma_1)$, $\mathcal P(\Gamma_2)$, one of which does not yield an absolute minimizer of $E$(in $\mathcal S(\Gamma_1)$ or $\mathcal S(\Gamma_2)$).     
\end{theorem}
{\bf Proof.} Indicate $\mf^i := \mf^i(x^i)$ for some $x^i\in M^{i\ast}, \,i=1,2$. For some $y\in \mm^{\ast}$, $\mf(y)$ is the strict relative minimum of $E$ in  $\mathcal S(\Gamma_1, \Gamma_2)$. Clearly, $y$ is also a strict relative minimizer of $\me$ in $\mm^{\ast}$. For $x=(x^1, x^2, 0)$, consider 
           \[ P = \{ p\in C^0([0,1], \overline{\mm}) | p(0) = x, p(1) = y \},\]
and         
            \[\beta := \inf_{p\in P}\max_{t\in[0,1]}\me(p(t)).\]
The PS condition implies that if $\beta > \max\{\me(x), \me(y)\}$, $\beta$ is a critical value which possesses a non-relative minimum critical point. Actually $\beta > \me(y)$, since $y$ is a strict relative minimizer. See \cite {s1} chapter II and \cite{ki} for details on that.\medskip

Supposing that any solution of $\mathcal P(\Gamma_i)$ is an absolute minimum of $E$ in $\mathcal S(\Gamma_i)$, we have a solution of $\mathcal P(\Gamma_1,\Gamma_2)$ which is not a relative minimum of $E$ in $\mathcal S(\Gamma_1,\Gamma_2)$, by the $E$-minimality of harmonic extensions.\medskip

It remains to show that $\beta := \inf_{p\in P}\max_{t\in[0,1]}\me(p(t)) > \me(x)$. We only need to consider $q=(q^1,q^2,\rho)\in p([0,1])$ for some $p\in P$ such that $\me(q^1,q^2,0) \le C$, $C$ a constant dependent on $N$.\medskip
    
Let $\varepsilon, \rho_0$ be as in Lemma \ref{boundary}, and consider the set of $q$'s with $\|q^i-{\tilde x}^i\|\ge \varepsilon$ for any absolute minimizer $\tilde x =({\tilde x}^1, {\tilde x}^2, 0)$ of $\me$ in $\partial \mm$. Then there exists $\delta_1>0$, dependent on $\varepsilon$, such that $\me(q^1,q^2,0)\ge \me(x)+\delta_1$ for all but finitely many $q$'s. If not, we would have a minimizing sequence converging to some absolute minimizer $\tilde x$ by the PS condition (Proposition \ref{psc}) and Proposition \ref{h22}, contradicting the choice of $q$.\medskip

Moreover, from the uniform convergence of $\me$  on a bounded set of $q^i$ (see Lemma \ref{differ}) when $\rho\rightarrow 0$, we can choose $\delta_2, \rho_1$ with $\delta_1-\delta_2> 0$, such that for all $\rho\in (0,\rho_1)$, $|\me(q^1,q^2,\rho)-\me(q^1,q^2,0)| \le \delta_2$.\medskip

Let $\bar \rho:=$min$\{\rho_0,\rho_1\}$. If $\|q^i-{\tilde x}^i\| < \varepsilon$ for some $ {\tilde x}$ as above, it is easy to see that $\me(q^1,q^2,\bar \rho) \ge \me(x) + \delta_3$ with $\delta_3>0$, by Lemma \ref{boundary}. If that were not so in fact, then $\me(q^1,q^2,\bar \rho)\ge \me(q^1,q^2,0) - \delta_2 \ge \me(x) + \delta_1-\delta_2$, by the above choices. This completes the proof.     \hfill $\Box$ \medskip
 
Now we specialize the main result to the three-dimensional sphere $S^3$ and hyperbolic space $H^3$, to which we can apply condition (C1).
\begin{example}
Let $\Gamma_1, \Gamma_2\subset B(p,\pi/2)$ for some $p\in S^3$, in other words $\Gamma_1, \Gamma_2$ lie in a hemisphere. Then the conclusion of the main theorem, under those conditions holds.
\end{example}
     
If there is exactly one solution to $\mathcal P(\Gamma_i)$, $i=1,2$, our main theorem guarantees that the existence of a minimal surface of annulus type whose energy is a strict relative minimum of $E$ in $\mathcal S(\Gamma_1, \Gamma_2)$ ensures the existence of an unstable minimal surface of annulus type. From \cite{lj}, a solution to $\mathcal P(\Gamma_i)$ is unique in $H\dach 3$ if the total curvature of $\Gamma_i$ is less than $4\pi$. Since $i(p)=\infty$ for all $p\in H^3$ we conclude 
\begin{example}
Let $\Gamma_1, \Gamma_2$ possess total curvature $\le 4\pi$ in $H^3$ and $\dist(\mf^1, \mf^2)>0$. If $E$ has a strict relative minimizer in $\mathcal S(\Gamma_1, \Gamma_2)$, then there exists an unstable minimal surface of annulus type in $H^3$.
\end{example}

\begin{appendix}

\section{Regularity of the critical points of $\me$}

This appendix is devoted to the proof of the following result, full details of which are found in \cite{ki2}.

\begin{theorem}\label{appendixthm}
Let $x=(x^1,x^2,\rho)\in M^1\times M^2 \times (0,1)$ with $g_i(x)=0,\, i=1,2$. Then $\mfr:=\mfr(x^1,x^2)$ belongs to $H^{2,2}(A_{\rho},N)$.  
\end{theorem}

Noting that $\mfr$ is harmonic in $N\stackrel{\eta}{\hookrightarrow}\mathbb R^k$, i.e. $\tau_h(f)=0$, polar coordinates give 
    \bqrs
        \lefteqn{ |\nabla^{2}\mfr|^2 =  |\partial_{r}d\mfr|^2 
                    + \frac{1}{r^2}|\partial_{\theta}d\mfr|^2}\hspace{0.5cm}\\
        & \le &  C(\varepsilon)|\Delta_{\mathbb R^k} \mfr|^2     
                + (2+\varepsilon)\frac{1}{r^2}|\partial_{\theta}d\mfr|^2
                +  C(\varepsilon)\frac{1}{r^2}\frac{1}{r^2}|\partial_{\theta}\mfr|^2\\
        & \le & C(\varepsilon,\eta,A_{\rho})|d\mfr|^2 + 
                   C(\varepsilon,\rho)|\partial_{\theta}d\mfr|^2.
    \eqrs
By a well known result of \cite{gt} it suffices to show that 
     \begin{equation} \label{endlich}
        \int_{A_\rho}|\Delta_{h}d\mfr|^2d\omega \le C < \infty, 
     \end{equation}
where $\Delta_{h}d\mfr:= \frac{d\mfr(r,\theta+h)-d\mfr(r,\theta)}{h}, h\not= 0$ and $C$ is independent of $h$. \medskip

Following \cite{ho}, observe that
\begin{remark}\label{hope}
For $\phi=(\phi^1,\phi^2)\in\hs\times\hs(\frac{\cdot}{\rho})$ define  
       \bqr
         \label{fifty}
        \mathbf A(\mfr)(\phi) := -\int_{A_{\rho}}\langle II\circ\mfr(d\mfr\,,d\mfr),X\rangle d\omega 
                                                    + \int_{A_{\rho}}\langle d\mfr,dX\rangle d\omega, 
     \eqr
where $X$ is any mapping in $H^{1,2}(A_{\rho},\mathbb R^k)$ with $X|_{\partial A_{\rho}} = \phi$. Then the expression on the right-hand-side only depends on the boundary of $X$, since $\mfr$ in harmonic in $N$.                \hfill $\Box$
\end{remark}

In particular, taking $\phi^{i}\in H^{\frac{1}{2},2}\cap C^{0}(\partial B,(x^{i})^{\ast}T\Gamma_{i}), i=1,2$, we consider $X:=\mj_{\mfr}(\phi^1,\phi^2)$, which is tangent to $N$ along $\mfr$. Since $\langle II\circ\mfr(d\mfr\,,d\mfr),\mj_{\rho}(\phi^1,\phi^2)\rangle\equiv 0$, 
   \bqr
     \mathbf A(\mfr)(\phi) 
        & = & \int_{A_{\rho}}\langle d\mfr, d\mj_{\mfr}(\phi^1,0)\rangle d\omega  
                 + \int_{A_{\rho}}\langle d\mfr, d\mj_{\mfr}(0,\phi^2)\rangle d\omega 
                 \label{div2} \\ 
        & = & \langle\partial_{x^1}\me,\phi^1\rangle + \langle\partial_{x^2}\me,\phi^2\rangle.
            \nonumber 
   \eqr
Hence for a critical point $x=(x^1,x^2,\rho)$ of $\me$,  $\mathbf A(\mfr)(\xi)\ge 0$ for all $\xi = (\xi^1,\xi^2)\in \mathcal T_{x^1}\times \mathcal T_{x^2}$.

\begin{lemma}\label{growthcondition}
For each $P_0\in \partial A_{\rho}$ there exist $C_0, \mu, r_0 >0$ such that for all $r\in [0,r_0]$ 
          \begin{equation}\label{sixteen}
               \int_{A_{\rho}\cap B_r(P_0)} (|d\mfr|^2 + |d\hr(\tilde{w^1},0)|^2)d\omega \le 
               C_0 r^{\mu}\int_{A_{\rho}} (|d\mfr|^2 + |d\hr(\tilde{w^1},0)|^2)d\omega. 
         \end{equation}
\end{lemma}
 
{\bf Proof of Lemma \ref{growthcondition}} \medskip

Let $P_0\in C_1$ fixed, define $B_r := B_r(P_0)$, and 
     \begin{equation*}
         \tilde{w^1_0}  :=      Q^{-1} 
          \int_{(B_{2r}\backslash B_r )\cap \partial B}\tilde{w^1} d_o,
          \quad w^1_0:= \tilde{w^1_0} + Id : \mathbb R \rightarrow \mathbb R, 
     \end{equation*}
where $ \int_{(B_{2r}\backslash B_r )\cap\partial B}d_o      := Q $. Then 
      \[\tilde{\xi}_{\phi} := -\big[ \phi(|e^{i\theta}-P_0|) \big]^2(w^1-w^1_0)
               \frac{\partial}{\partial\theta}\circ \bar{w^1} \in 
                \hs\cap C^0 (\partial B, \bar{w^1}^{\ast}T(\partial B)),\]
where $\bar{w^1}$ is a map from $\partial B$ to itself, and $\phi\in C^{\infty}$ is a non-increasing function of $|z|$ satisfying $ 0\le \phi(z)\le 1$, $\phi\equiv 1$ if $|z| \le 2r$, $\phi\equiv 0$ for $|z|\ge 3r$, plus $|d \phi|\le\frac{C}{r}$ and $|d^2 \phi|\le \frac{C}{r^2}$.\medskip 

Since $(1-\phi^2)w^1+\phi^2 w^1_0\in W^1_{\mathbb R^k}$, we see that $d\gamma^1(\tilde{\xi}_{\phi}) \in \mathcal T_{x^1}$, and $\mathbf A(\mfr)(d\gamma^1(\tilde{\xi}_{\phi}),0)\ge 0$. \medskip
      
For $x^1_0 := \gamma^1(w^1_0)$, 
    \bqrs
       x^1-x^1_0 & = & d\gamma^1(w^1-w^1_0) - 
           \underbrace{\int^{w^1}_{w^1_0}\int^{w^1}_{s'}d^2\gamma^1(s'')ds''ds'}_{=:\alpha(w^1)}, 
    \eqrs
and for small $r>0$,
     \bqrs
      \mathbf A(\mfr)(\phi^2(\mfr-\mfr^0)|_{C_1}, 0) 
          & = & \mathbf A(\mfr)(\phi^2 d\gamma^1(w^1-w^1_0),0) 
                  - \mathbf A(\mfr)(\phi^2 \alpha(w^1),0)\\
          & \le & -\mathbf A(\mfr)(\phi^2 \alpha(w^1),0),
     \eqrs
where $\mfr^0(A_{\rho})\equiv x^1_0\in \Gamma_1$. \medskip

On the other hand, for small $r>0, \phi^2(\mfr-\mfr^0)|_{C_2}\equiv 0$, so we can take $\phi^2(\mfr-\mfr^0)$ instead of $\phi$ in the definition of $\mathbf A(\mfr)$, to the effect that  
   \bqr
     \lefteqn{\int_{A_{\rho}}\langle \phi^2d\mfr,d\mfr\rangle d\omega \le   
          \int_ {A_{\rho}}\langle \phi^2(\mfr-\mfr^0), 
                       II\circ \mfr(d\mfr,d\mfr)\rangle d\omega}\nonumber\\ 
          & & \hspace{4.0cm} - \int_ {A_{\rho}}\langle 2\phi d\phi(\mfr-\mfr^0),d\mfr\rangle d\omega 
                        -\mathbf A(\mfr)(\phi^2 \alpha(w^1),0).\label{nine}
    \eqr 

Now define a real valued map of $(r,\theta)\in [\rho,1]\times \mathbb R$ as follows: 
       \[ T^1(w^1)(r,\theta) := H_{\rho}(\tilde w, 0)(r,\theta) + Id(r,\theta) \en \text{with} 
          \en Id(r,\theta) = \theta, \]
where $H_{\rho}(\tilde w, 0)$ is the harmonic extension to $A_{\rho} \approx [\rho,1]\times \mathbb R/2\pi$ with $\tilde w$ on  $\partial B$ and $0$ on $\partial B_{\rho}$.\medskip

In order to estimate $-\mathbf A(\mfr)(\phi^2 \alpha(w^1),0)$, we consider 
       \[\widetilde{\star\star}
         :=\phi^2\int^{T^1(w^1)}_{w^1_0}\int^{T^1(w^1)}_{s'}d^2\gamma^1(s'')ds''ds'\in 
         H^{1,2}(A_{\rho},\mathbb R^k)\]
with $\widetilde{\star\star}|_{C_1}= \phi^2 \alpha(w^1), \widetilde{\star\star}|_{C_2}\equiv 0$, where $w^1_0(r,\theta) = \tilde{w^1_0} + Id(r,\theta) = \tilde{w^1_0} +\theta, \en (r,\theta)\in [\rho,1]\times \mathbb R$.\medskip

an easy computation shows that 
      \bqrs
        |\widetilde{\star\star}| 
          & \le & C(\gamma^1,x^1) \phi^2 |\hr(\tilde{w^1},0)-\tilde{w^1_0}|^2, \\
                                                   |d \widetilde{\star\star}| 
          & \le & C(\gamma^1,x^1) |\hr(\tilde{w^1},0)-\tilde{w^1_0}|^2 \phi|d\phi| 
                         + C(\gamma^1,x^1) |d\hr(\tilde{w^1},0)||\hr(\tilde{w^1},0)-\tilde{w^1_0}|^2 \phi^2,  
      \eqrs
and from (\ref{nine}) Young's inequality implies 
   \bqrs
        \lefteqn{\int_{A_{\rho}}\langle \phi^2d\mfr,d\mfr\rangle d\omega \le  \int_  
                           {A_{\rho}}|d\mfr|^2|\mfr-\mfr^0|\phi^2d\omega }\\
              && \hspace{2.0cm}+ \frac{\varepsilon}{5}\int_{A_{\rho}}|d\mfr|^2 \phi^2 d\omega 
                            + C(\varepsilon)\int_{A_{\rho}}|\mfr-\mfr^0|^2|d\phi|^2  d\omega \\
              && \hspace{2.0cm}+  C\|\hr(\tilde{w^1},0)-\tilde{w^1_0}\|_{L^{\infty}(B_{3r})} 
              \int_{A_{\rho}}\big( |d\mfr|^2 \phi^2 
                       +  |\hr(\tilde{w^1},0)-\tilde{w^1_0}|^2 |d\phi|^2 \big)d\omega\\
              && \hspace{2.0cm}+  C\|\hr(\tilde{w^1},0)-\tilde{w^1_0}\|_{L^{\infty}(B_{3r})}
                    \int_{A_{\rho}}\big( |d \hr(\tilde{w^1},0)|^2 + |d\mfr|^2 \big) \phi^2 d\omega \\
              && \hspace{2.0cm}+  C \int_{A_{\rho}} |\hr(\tilde{w^1},0)
                  -\tilde{w^1_0}|^2 |d\mfr|^2 \phi^2 d\omega.   
     \eqrs   
Thus, for a sufficiently small $r\in (0, r_0)$ dependent on $\varepsilon$, $C$, and the modulus of continuity of $\mfr-\mfr^0$ and $\hr(\tilde{w^1},0)-\tilde{w^1_0} $, we have an estimate:
      \bqr   
         \lefteqn{\int_{A_{\rho}}\langle \phi^2d\mfr,d\mfr\rangle d\omega \le 
           \varepsilon\int_{A_{\rho}} \big( |d\mfr|^2  
                  +  |d\hr(\tilde{w^1},0)|^2 \big)\phi^2 d\omega}\nonumber \\
           & & \hspace{4.0cm} +C(\varepsilon)\int_{A_{\rho}} \big( |\mfr-\mfr^0|^2 
                  +  |\hr(\tilde{w^1},0)-\tilde{w^1_0}|^2 \big)|d\phi|^2 d\omega.\label{ten}
       \eqr

This corresponds to (5.6) in \cite{s1} (Proposition 5.1, II). A completely similar computation for $\int_{A_{\rho}} |d\hr(\tilde{w^1},0)|^2 \phi^2 d\omega$ and $\int_{ A_{\rho}\cap B_r(P_0)} \big( |d\mfr|^2+|\hr(\tilde{w^1},0)|^2 \big)d\omega$ eventually yields (\ref{sixteen}). \hfill $\Box$ \bigskip
 
{\bf Proof of Theorem \ref{appendixthm}} \medskip

We will show (\ref{endlich}) by several steps.\medskip 

{\bf (I)} With $\Delta_{-h}\Delta_{h}\mfr|_{\partial B} = \Delta_{-h}\Delta_{h}\gamma^1\circ e^{iw^{1}}$ and $\Delta_{-h}\Delta_{h}\mfr|_{\partial B_{\rho}}(\cdot\rho) = \Delta_{-h}\Delta_{h}\gamma^2\circ e^{iw^{2}(\cdot)}$,
      \bqrs
        \lefteqn{\int_{A_{\rho}}|\Delta_hd\mfr|^2d\omega        = 
           -\int_{A_{\rho}}\langle d\mfr,d\Delta_{-h}\Delta_{h}\mfr\rangle d\omega}\hspace{1.0cm}\\
           & = & -\int_{A_{\rho}}\langle II\circ \mfr(d\mfr,d\mfr),
                  \Delta_{-h}\Delta_h\mfr\rangle d\omega
                  -\mathbf A(\mfr)(\Delta_{-h}\Delta_{h}\mfr|_{\partial A_{\rho}}).
     \eqrs                                                                                        
         
Denoting $\gamma^1\circ e^{iw^1}$, $\gamma^2\circ e^{iw^2}$ by $\gamma^1(w^1(\theta))$, $\gamma^2(w^2(\theta))$ and $w^i(\cdot + h)$, $w^i(\cdot - h)$ by $w^i_{+}$, $w^i_{-}$ respectively, we have:
     \bqrs
       \lefteqn{\Delta_{-h}\Delta_h\gamma^i(w^i) 
               =  \Delta_{-h}\left[d\gamma^i(w^i)\left(\frac{w^i_{+}-w^i_{-}}{h}\right)+
                 \frac{1}{h}\int^{w^i_{+}}_{w^i}\int^{s'}_{w^i}d^2\gamma^i(s'')ds''ds'\right]}\\
           & = & d\gamma^i(w^i)(\Delta_{-h}\Delta_{h}w^i)
               \underbrace{-\frac{1}{h}\int^{w^i_{-}}_{w^i}d^2\gamma^i(s')ds'\cdot\Delta_{h}w^i_{-}
               + \Delta_{-h}\left(\frac{1}{h}\int^{w^i_{+}}_{w^i}
                  \int^{s'}_{w^i}d^2\gamma^i(s'')ds''ds'\right)}_{=:P^i}
    \eqrs

Clearly $d\gamma^i(w^i)(\Delta_{-h}\Delta_{h}w^i)\in H^{\frac{1}{2},2} \cap C^{0}(\partial B,(x^{i})^{\ast}T\Gamma_{i})$.\medskip

Now define a map $S(P^1,0):A_{\rho}\rightarrow \mathbb R^k$ with boundary $(P^1,0)$ as follows: 
      \bqrs
         S(P^1,0)
            := -\frac{1}{h}\int^{T^1(w^1_{-})}_{T^1(w^1)(\cdot)} d^2\gamma^1(s')ds'\cdot 
                            H_{\rho}(\Delta_h w^1_{-},0) +\Delta_{-h}\big(\frac{1}{h} 
             \int^{T^1(w^1_{+})}_{T^1(w^1)}\int^{s'}_{T^1(w^1)} d^2\gamma^1(s'')ds''ds'\big).
      \eqrs
  
Similarly, we have a map $S(0,P^2)(\cdot):A_{\rho}\rightarrow \mathbb R^k$ with $0$ on $C_1$ and $P^2$ on $C_2$.\medskip   
        
By computation, $\frac{h^2}{2}\Delta_{-h}\Delta_{h}w^i =\frac{1}{2}(w^i_{-} +w^i_{+}) - w^i$, and $\frac{1}{2}(w^i_{-} +w^i_{+})\in W^i_{\mathbb R^k}$, noting that $W^i_{\mathbb R^k}$ is convex. Thus $\frac{h^2}{2}d\gamma^i(w^i)(\Delta_{-h}\Delta_{h}w^i) \in \mathcal T_{x^i}$ by definition of $\mathcal T_{x^i}$,. \medskip

As $g^1(x)=g^2(x)=0$, we have
        \bqr
           \mathbf A(\mfr)\left( d\gamma^1(w^1)(\Delta_{-h}\Delta_{h}w^i),
                     d\gamma^2(w^2)(\Delta_{-h}\Delta_{h}w^2)\right)
          \ge 0, \label{sixty}
        \eqr
where we have dropped the scaling term $(\frac{\cdot}{\rho})$ in relative to second variation. Now  
      \bqr         
         \lefteqn{\int_{A_{\rho}}|\Delta_hd\mfr|^2d\omega 
              = -\int_{A_{\rho}}\langle II\circ \mfr(d\mfr,d\mfr),
                     \Delta_{-h}\Delta_h\mfr\rangle d\omega
                -\mathbf A(\mfr)(\Delta_{-h}\Delta_{h}\mfr|_{\partial A_{\rho}})}
                    \hspace{1.0cm}\nonumber\\
          &=& -\int_{A_{\rho}}\langle II\circ \mfr(d\mfr,d\mfr),
             \Delta_{-h}\Delta_h\mfr\rangle d\omega \nonumber\\
          & & \hspace{1.0cm} -\mathbf A(\mfr)(P^1,P^2)
                  -\mathbf A(\mfr)\left(d\gamma^1(w^1)(\Delta_{-h}\Delta_{h}w^1),
                d\gamma^2(w^2)(\Delta_{-h}\Delta_{h}w^2)\right)\nonumber\\
        & \le & -\int_{A_{\rho}}\langle II\circ \mfr(d\mfr,d\mfr),
                \Delta_{-h}\Delta_h\mfr\rangle d\omega -\mathbf A(\mfr)(P^1,P^2)\nonumber\\
        & \le & -\int_{A_{\rho}}\langle II\circ \mfr(d\mfr,d\mfr),
                \Delta_{-h}\Delta_h\mfr\rangle d\omega \label{first}\\
        & & + \int_{A_{\rho}}\langle II\circ \mfr(d\mfr,d\mfr),S(P^1,0)\rangle d\omega 
            + \int_{A_{\rho}}\langle II\circ \mfr(d\mfr,d\mfr),S(0,P^2)\rangle d\omega
                      \label{second}\\
        & & - \int_{A_{\rho}}\langle d\mfr,dS(P^1,0)\rangle d\omega 
            - \int_{A_{\rho}}\langle d\mfr,dS(0,P^2)\rangle d\omega.\label{third}
        \eqr
 
For the estimates of these terms we need some preliminaries.\medskip
        
First, let $s(\tau):=\tau\mf_{\rho,+}+(1-\tau)\mfr,\, 0\le \tau\le 1$. Then 
    \bqrs
     \lefteqn{|\Delta_h II\circ \mfr (d\mfr,d\mfr)| 
        = |\frac{1}{h}\{ II\circ \mf_{\rho,+}(\mf_{\rho,+}, \mf_{\rho,+}) 
                        - II\circ \mfr(d\mfr,d\mfr)\}|}\\
          &=& |\frac{1}{h}\{II\circ\mf_{\rho,+}(d\mf_{\rho,+},d\mf_{\rho,+})
                   -II\circ\mfr(d\mf_{\rho,+},d\mf_{\rho,+})
         + II\circ\mfr(d\mf_{\rho,+},d\mf_{\rho,+}) 
                   - II\circ\mfr(d\mfr,d\mfr)\}|\hspace{2.0cm}\\
          &=& |dII(\mfr)\cdot\Delta_h \mfr(d\mf_{\rho,+},d\mf_{\rho,+})
                    + \frac{1}{h}\int^{1}_{0}\int^{t}_{0}d^2II(s(\tau))
                     |\mf_{\rho,+}-\mfr|^2d\tau dt(d\mf_{\rho,+},d\mf_{\rho,+})\\
            & &\hspace{2.0cm} + II\circ\mfr(\Delta_h d\mfr,d\mf_{\rho,+}) 
                  + II\circ\mfr(d\mfr,\Delta_h d\mfr)|\\
           &\le& C(\|\mfr\|_{C^0(\ar)})[|\Delta_h \mfr||d\mf_{\rho,+}|^2 
                  + |\Delta_h d\mfr|(|d\mf_{\rho,+}| + |d\mfr|)].
      \eqrs 
   
Letting 
        \[-\frac{1}{h}\int^{T^1(w^1_{-})}_{T^1(w^1)}d^2\gamma^1(s')ds'
                := \star \quad  \text{and}  \quad
          \frac{1}{h}\int^{T^1(w^1_{-})}_{T^1(w^1)}\int^{s'}_{T^1(w^1)}d^2\gamma^1(s'')ds''ds'
           :=\star\star,\]
we have                              
        \[ |\star| \le C(\gamma^1)|H_{\rho}(\Delta_{-h}w^1,0)|,  
               \quad |\star\star|\le C(\gamma^1)|H_{\rho}(\Delta_{h}w^1,0)|,\]
and
        \bqrs
           |d\star| 
            & \le & C(\|\gamma^1\|_{C^3})\big(|H_{\rho}(\Delta_{-h}w^1,0)||dH_{\rho}(w^1_{-},0)| 
                   + |dH_{\rho}(\Delta_{-h}w^1,0)|\big), \\
           |d\star\star|  
                     & \le & C(\|\gamma^1\|_{C^2})|H_{\rho}(\Delta_{h}w^1,0)|
                             \big(|dH_{\rho}(\tilde{w^1}_{+},0)| + |dH_{\rho}(\tilde{w^1},0)|\big). 
       \eqrs   

With all that, we can estimate (\ref{first}), (\ref{second}), (\ref{third}) using a $C\in \mathbb R$, independent of $h$. All-in-all then, 
       \bqr
         \int_{A_{\rho}}|\Delta_h d\mfr|^2 d\omega
          &=& \varepsilon C\int_{A_{\rho}}|\Delta_h d\mfr|^2d\omega
                + \varepsilon C\int_{A_{\rho}}|d\hr(\Delta_h w^1,0)|^2d\omega 
                  +C(\varepsilon)\Xi \, \label{fourth},
       \eqr

where $\Xi$ stands for: 
 \bqrs
    \lefteqn {
        \int_{A_{\rho}}\big( 
               |d\hr(\tilde{w}^1_{-},0)|^2+ |d\hr(\tilde{w}^1_{+},0)|^2
                  +|d\hr(\tilde{w}^1,0)|^2 +|d\hr(0,\tilde{w}^2_{-})|^2+ 
                    |d\hr(0,\tilde{w}^2_{+})|^2 + d\hr(0,\tilde{w}^2)|^2  } \\
         & & +||d\mfr|^2  \big) 
             \cdot \big(|\Delta_h\mfr|^2+|\hr(\Delta_{-h}w^1,0)|^2 
                + |\hr(\Delta_{h}w^1,0)|^2
                +|\hr(0,\Delta_{-h}w^2)|^2+ |\hr(0,\Delta_{h}w^2)|^2\big)d\omega.
 \eqrs

{\bf (II)} On  $\partial B$ we know that $\Delta_h(\gamma^i\circ w^i) = d\gamma^i(w^i)\Delta_hw^i + \frac{1}{h}\int^{w^i_{+}}_{w^i}\int^{s'}_{w^i}d^2\gamma^i(s'')ds''ds' $, so  
          \begin{equation}\label{fifth} 
               \Delta_h w^i = |d\gamma^i(w^i)|^{-2}\big[ d\gamma^i(w^i)\cdot \Delta_h\mfr 
                - d\gamma^i(w^i)\cdot\frac{1}{h}\int^{w^i_{+}}_{w^i}\int^{s'}_{w^i}d^2
                          \gamma^i(s'')ds''ds'\big]. 
          \end{equation}

Using $T^i(w^i)$ on the right-hand-side of (\ref{fifth}), we get an $H^{1,2}(A_{\rho},\mathbb R^k)$-extension with boundary $\Delta_h w^i$ on $C^1$ and $0$ on $C_2$, and by D-minimality of the harmonic extension among the maps with same boundary, it follows that 
        \bqr
           \int_{A{\rho}}|d\hr(\Delta_h w^1,0)|^2 d\omega
            & \le & C\int_{A_{\rho}} \big[ |d\hr(w^1,0)| \big( |\Delta _h\mfr|+|\star\star| \big) + 
                           |d\Delta_h \mfr| + |d\star\star| \big]^2 d\omega \nonumber\\
            & \le &  C\int_{A_{\rho}}|d\Delta_h \mfr|^2 d\omega  + C\Xi,\,\label{sixth}
        \eqr
again by Young's inequality. We van attain a similar estimate for $\int_{A{\rho}}|d\hr(0,\Delta_h w^2)|^2 d\omega$.\medskip

Using the estimate (\ref{fourth}) for $\int_{A_{\rho}}|d\Delta_h \mfr|^2 d\omega$, (\ref{sixth}) implies  
      \bqrs
           \lefteqn{\int_{A_{\rho}}|d\Delta_h \mfr|^2 d\omega 
                + \int_{A_{\rho}}|d\hr(\Delta_h w^1,0)|^2 d\omega
                + \int_{A_{\rho}}|d\hr(0,\Delta_h w^2)|^2 d\omega}\\
            & \le & \varepsilon C \int_{A_{\rho}}|d\Delta_h \mfr|^2 d\omega 
                      + \varepsilon C\int_{A_{\rho}}|d\hr(\Delta_h w^1,0)|^2 d\omega
                      + \varepsilon C\int_{A_{\rho}}|d\hr(0,\Delta_h w^2)|^2 d\omega 
                       + C(\varepsilon) \Xi\,.
      \eqrs

For some small $\varepsilon>0$ in the above formula we get the inequality:
       \bqrs
        \lefteqn{\int_{A_{\rho}}|\Delta_h d\mfr|^2 d\omega 
                       + \int_{A_{\rho}}|d\hr(\Delta_h w^1,\Delta_h w^2)|^2 d\omega}\hspace{1.0cm}\\
                    & \le & C(\varepsilon) \int_{A_{\rho}} 
                       \big( |d\mfr|^2 + |d\mf_{\rho +}|^2 + |d\mf_{\rho -}|^2 
                          + |d\hr(\tilde{w^1},\tilde{w^2}|)^2 
                          + |d\hr(\tilde{w^1_+},\tilde{w^2_+})|^2 \nonumber\\
                    & & \hspace*{1.0cm}  
                         + |d\hr(\tilde{w^1_{-}},\tilde{w^2_-}|^2 \big)\cdot 
                          \big( |\Delta_h \mfr|^2 + |H(\Delta_{-h}w^1,\Delta_{-h}w^2)|^2 + 
                                   |H(\Delta_{h}w^1,\Delta_{h}w^2)|^2 \big) d\omega \nonumber.
        \eqrs

Extend now $\mfr$ to $\mathbb R^2\backslash B_{\rho^2}$ by conformal reflection 
        \bqrs
           \mfr(z) & = & \mfr\big( \frac{z}{|z|^2} \big) ,
                         \en \text{if} \en 1\le |z| \\
           \mfr(z) & = & \mfr\big( \frac{z}{|z|^2}\rho^2 \big) ,
                         \en \text{if} \en \rho^2 \le |z|\le \rho.
       \eqrs 

Choose $r\in \big( 0, \min\{\frac{\rho-\rho^2}{2}, r_0\} \big)$, and $ \varphi \in C^{\infty}_0\big(B_{2r}(0) \big)$ with $\varphi \equiv 1$ on $B_{r}(0)$.\medskip

We may cover $A_{\rho}$ with balls of radius $r$ in such a way that any $p\in A_{\rho}$ lies in the intersection of at most $k$ balls, for any $r$ as above (recall $\mathbb R^2$ is metrizable). Let $B^i$ denote the balls with centres $p_i$ and set $\varphi_i(p):= \varphi(p-p_i)$. Then 
     \bqrs                        
        \lefteqn{\int_{A_{\rho}}|\Delta_h d\mfr|^2 d\omega 
                    + \int_{A_{\rho}}|d\hr(\Delta_h w^1,\Delta_h w^2)|^2 d\omega}\\
        & \le & C\,\Sigma_{i} \int_{\mathbb R^2\backslash A_{\rho^2}}
                   \big( |\Delta_h \mfr|^2 + |H(\Delta_{-h}w^1,\Delta_{-h}w^2)|^2 + 
                   |H(\Delta_{h}w^1,\Delta_{h}w^2)|^2 \big)\varphi_i^2 \cdot \\
        & &   \underbrace{\big( |d\mfr|^2 + |d\mf_{\rho +}|^2   + |d\mf_{\rho -}|^2 
                + |d\hr(\tilde{w^1},\tilde{w^2})|^2 + |d\hr(\tilde{w^1_+},\tilde{w^2_+})|^2
                + |d\hr(\tilde{w^1_{-}},\tilde{w^2_-})|^2 \big)}_{=: \chi}d\omega.
      \eqrs

By substituting $|d\mf_{\rho +}|^2$ and $|d\hr(\tilde{w^1_+},\tilde{w^2_+})|^2$ (or $|d\mf_{\rho -}|^2$ and $|d\hr(\tilde{w^1_-},\tilde{w^2_-})|^2$) in Lemma \ref{growthcondition}, we conclude that $\chi$ satisfies the growth condition of Morrey. Now applying Morrey's Lemma (1, Lemma 5.4.1 \cite{mo}) to $\chi$ and $ (\Delta_h \mfr)\varphi_i $, $\chi$ and $H(\Delta_{-h}w^1,\Delta_{-h}w^2)\varphi_i$ or $\chi$ and $H(\Delta_{h}w^1,\Delta_{h}w^2)\varphi_i$, and adding over the index $i$ for some small $r>0$, we obtain a constant $C>0$, independent of $|h|\le h_0$, such that      
         \[ \int_{A_{\rho}}|\Delta_h d\mfr|^2 d\omega  \le C. \] 
                                                                        \hfill $\Box$

\end{appendix}


\begin{thebibliography}{30}
\bibitem[Ki2]{ki2} Kim, H., A variational approach to the regularity of minimal surfaces of annulus type in Riemmanian manifolds, To appear in Differential Geometry and its Application
\end{thebibliography}
\end{document}